# CONTROL LYAPUNOV FUNCTIONS AND STABILIZATION BY MEANS OF CONTINUOUS TIME-VARYING FEEDBACK

**Iasson Karafyllis**[*] **and John Tsinias**[**]

[*]**Department of Environmental Engineering, Technical University of Crete, 73100, Chania, Greece**
email: ikarafyl@enveng.tuc.gr

[**]**Department of Mathematics, National Technical University of Athens, Zografou Campus 15780, Athens, Greece**
Email: jtsin@central.ntua.gr

**Abstract**

For a general time-varying system, we prove that existence of an "Output Robust Control Lyapunov Function" implies existence of continuous time-varying feedback stabilizer, which guarantees output asymptotic stability with respect to the resulting closed-loop system. The main results of the present work constitute generalizations of a well known result towards feedback stabilization due to J. M. Coron and L. Rosier concerning stabilization of autonomous systems by means of time-varying periodic feedback.

## 1. Introduction

Control Lyapunov functions play a central role to the solvability of the feedback stabilization problem and several important works are found in the literature, where sufficient conditions are provided in terms of Lyapunov functions for characterizations of various types of stability, as well for existence of feedback stabilizers (see for instance [1,2,4,5,6,7,9,11,12,13,15,18,19,20,21,25]). In the present work we consider time-varying uncertain systems of the general form:

$$\dot{x} = f(t, d, x, u)$$
$$Y = H(t, x) \qquad (1.1)$$
$$x \in \Re^n, d \in D, t \geq 0, Y \in \Re^k, u \in U \subseteq \Re^m$$

where $d(\cdot)$ is a time-varying disturbance which takes values on the set $D \subset \Re^l$ and $u, Y$ play the role of input and output, respectively, of (1.1). We assume that $0 \in \Re^n$ is an equilibrium for (1.1), i.e. it holds $f(t,d,0,0) = 0$ and $H(t,0) = 0$ for all $(t,d) \in \Re^+ \times D$. We prove that existence of an "Output Robust Control Lyapunov Function" (ORCLF) implies existence of continuous time-varying feedback stabilizer

$$u = K(t, x) \qquad (1.2)$$

that guarantees global output asymptotic stability of the output $Y = H(t,x)$ with respect to the resulting closed-loop system (1.1) with (1.2), being uniform with respect to disturbances $d(\cdot)$. Our main results constitute generalizations of an important result towards feedback stabilization obtained in [7] by J. M. Coron and L. Rosier concerning autonomous systems:

$$\dot{x} = f(x, u), \ (x, u) \in \Re^n \times \Re^m \text{ with } f(0,0) = 0 \in \Re^n \qquad (1.3)$$

Particularly, among other things in [7], it is established that existence of a time-independent control Lyapunov function, which satisfies the "small-control property" guarantees existence of a continuous time-varying periodic feedback (1.2) in such a way that $0 \in \Re^n$ is globally asymptotically stable for the resulting time-varying closed-loop system (1.1) with (1.2). In the present work we present generalizations of the result above for general time-varying systems (1.1). Particularly, in Theorem 2.8 of present paper we establish that existence of an ORCLF, which satisfies a time-varying version of the small control property, implies existence of a continuous feedback stabilizer $K : \Re^+ \times \Re^n \to U$, which is continuously differentiable with respect to $x \in \Re^n$ on the set $\Re^+ \times (\Re^n \setminus \{0\})$ exhibiting Robust Global Asymptotic Output Stability (RGAOS) of the resulting closed-loop system, being *uniform with respect to initial values of time*. In Theorem 2.9 of present work it is shown that, under lack of the small control property, existence of an ORCLF, implies existence of a continuous feedback stabilizer $K : \Re^+ \times \Re^n \to U$, being continuously differentiable with respect to $x \in \Re^n$ on the set $\Re^+ \times \Re^n$ exhibiting RGAOS of the resulting closed-loop system, being in general *non-uniform with respect to initial values of time*. We note here, that various concepts of asymptotic stability being in general non-uniform with respect to initial values of time, their Lyapunov characterizations, as well applications to feedback stabilization and related problems are found in several recent works (see for instance [11,12,13] and references therein). As a consequence of Theorem 2.9 and the main result in [12], it is shown in Corollary 2.10 of the present work that that the converse claim of Theorem 2.9 is true; to be more precise the following three statements are equivalent:

• *existence of an ORCLF (under lack of the small control property),*

• *existence of a continuous mapping* $K : \Re^+ \times \Re^n \to U$ *being continuously differentiable with respect to* $x \in \Re^n$ *on* $\Re^+ \times \Re^n$, *such that the closed-loop system (1.1) with (1.2) is (non-uniformly in time) RGAOS,*

• *existence of an ORCLF satisfying the small control property.*

It should be emphasized here that, when the result of Theorem 2.9 is restricted to autonomous systems (1.3), we get the following result which generalizes both Artstein's theorem on stabilization in [2,20] and Rosier-Coron main result in [7]: Assume that (1.3) possess a (time-independent) Control Lyapunov Function, namely, suppose that there exists a map $V \in C^1(\Re^n; \Re^+)$, functions $a_1, a_2 \in K_\infty$, and $\rho \in C^0(\Re^+; \Re^+)$ being positive definite such that

$$a_1(|x|) \leq V(x) \leq a_2(|x|), \quad \min_{u \in U} \frac{\partial V}{\partial x}(x) f(x,u) \leq -\rho(V(x)), \quad x \in \Re^n.$$

Then there exists a time-varying continuous mapping $K : \Re^+ \times \Re^n \to U$, being continuously differentiable with respect to $x \in \Re^n$ on $\Re^+ \times \Re^n$, such that the closed-loop time-varying system (1.3) with (1.2) is RGAOS in general non-uniformly with respect to initial values of time.

Comparing with the results obtained in [7] the result above presents the important advantage that feedback stabilization is exhibited under lack of the small control property and the corresponding feedback is an ordinary map, being in general time-varying but non-periodic. On the other hand, our approach leads in general to non-uniform in time asymptotic stability for the resulting closed-loop system. Finally, it should be pointed out that the main results in the present work (Theorem 2.8 and Theorem 2.9) generalize the main result in [7] in the following additional directions:

- The dynamics of systems we consider are in general time-varying, including disturbances, and the control set $U$ is in general a positive cone of $\Re^m$.
- The general problem of robust output stabilization is considered and feedback stabilization is exhibited under the presence of time-varying Control Lyapunov Functions.

The proofs of the main results in the present work are inspired by the proof of the main result in [7], but are essentially different in many points.

The paper is organized as follows. In Section 2, several stability notions and the concept of the Output Robust Control Lyapunov Function are presented, as well precise statements of our main results are given. Section 3 contains the proofs of the main results.

**Notations** Throughout this paper we adopt the following notations:
∗ Let $A \subseteq \Re^n$. By $C^0(A; \Omega)$, we denote the class of continuous functions on $A$, which take values in $\Omega$. Likewise, $C^1(A; \Omega)$ denotes the class of functions on $A$ with continuous derivatives, which take values in $\Omega$.

* For a vector $x \in \Re^n$ we denote by $|x|$ its usual Euclidean norm and by $x'$ its transpose.
* $Z$ denotes the set of integers, $Z^+$ denotes the set of non-negative integers and $\Re^+$ denotes the set of non-negative real numbers.
* We denote by $[r]$ the integer part of the real number $r$, i.e., the greatest integer, which is less than or equal to $r$.
* A continuous mapping $\Re^+ \times \Re^n \ni (t,x) \to k(t,x) \in U$, is *continuously differentiable with respect to* $x \in \Re^n$ on the open set $A \subseteq \Re^+ \times \Re^n$ (with respect to the $\Re^+ \times \Re^n$ topology), if the mapping $A \ni (t,x) \to \frac{\partial k}{\partial x}(t,x) \in \Re^n$ is continuous and is called *locally Lipschitz with respect to* $x \in \Re^n$ on the open set $A \subseteq \Re^+ \times \Re^n$, if for every compact set $S \subseteq A$ it holds that

$$\sup \left\{ \frac{|k(t,x) - k(t,y)|}{|x-y|} : (t,x) \in S, (t,y) \in S, x \neq y \right\} < +\infty$$

* We denote by $K^+$ the class of positive $C^0$ functions defined on $\Re^+$. We say that a function $\rho : \Re^+ \to \Re^+$ is positive definite if $\rho(0) = 0$ and $\rho(s) > 0$ for all $s > 0$. By $K$ we denote the set of positive definite, increasing and continuous functions. We say that a positive definite, increasing and continuous function $\rho : \Re^+ \to \Re^+$ is of class $K_\infty$ if $\lim_{s \to +\infty} \rho(s) = +\infty$.
* Let $D \subseteq \Re^l$ be a non-empty set. By $M_D$ we denote the class of all Lebesgue measurable and locally essentially bounded mappings $d : \Re^+ \to D$.

## 2. Basic Notions and Main Results

In this work, we consider systems of the form (1.1) under the following hypotheses:

**(H1)** The mappings $f : \Re^+ \times D \times \Re^n \times U \to \Re^n$, $H : \Re^+ \times \Re^n \to \Re^k$ are continuous and for every bounded interval $I \subset \Re^+$ and every compact set $S \subset \Re^n \times U$ there exists $L \geq 0$ such that $|f(t,d,x,u) - f(t,d,y,v)| \leq L|x-y| + L|u-v|$ for all $(t,d) \in I \times D$, $(x,u) \in S$, $(y,v) \in S$.

**(H2)** The set $D \subset \Re^l$ is compact and $U$ is a closed positive cone, i.e., $U \subseteq \Re^m$ is a closed set and, if $u \in U$, then $(\lambda u) \in U$ for all $\lambda \in [0,1]$.

**(H3)** Zero $0 \in \Re^n$ is an equilibrium; particularly, assume that $f(t,d,0,0) = 0$, $H(t,0) = 0$ for all $(t,d) \in \Re^+ \times D$.

**Definition 2.1** *We say that a function $V \in C^1(\Re^+ \times \Re^n; \Re^+)$ is an **Output Robust Control Lyapunov Function (ORCLF) for (1.1)**, if there exist functions $a_1, a_2 \in K_\infty$, $\mu, \beta \in K^+$, $\rho \in C^0(\Re^+; \Re^+)$, being locally Lipschitz and positive definite and $b \in C^0(\Re^+ \times (\Re^n \setminus \{0\}); \Re^+)$, such that*

i) for every $(t,x) \in \Re^+ \times \Re^n$ it holds:

$$a_1(|H(t,x)| + \mu(t)|x|) \leq V(t,x) \leq a_2(\beta(t)|x|) \qquad (2.1)$$

ii) for every $(t,x) \in \Re^+ \times (\Re^n \setminus \{0\})$ it holds:

$$\min_{\substack{|u| \leq b(t,x) \\ u \in U}} \max_{d \in D} \left( \frac{\partial V}{\partial t}(t,x) + \frac{\partial V}{\partial x}(t,x) f(t,d,x,u) \right) \leq -\rho(V(t,x)) \qquad (2.2)$$

*For the case, where (2.2) holds and, instead of (2.1), it holds:*

$$a_1(|x|) \leq V(t,x) \leq a_2(\beta(t)|x|), \quad \forall (t,x) \in \Re^+ \times \Re^n \tag{2.1'}$$

the corresponding $V \in C^1(\Re^+ \times \Re^n; \Re^+)$ is called **State Robust Control Lyapunov Function (SRCLF).**

We say that the ORCLF (SRCLF) satisfies the **small-control property** with respect to (1.1), if in addition to (2.1), (2.2), ((2.1'), (2.2)), there exist functions $a_3 \in K_\infty$, $\gamma \in K^+$ such that the following inequality holds for all $(t,x) \in \Re^+ \times (\Re^n \setminus \{0\})$ :

$$b(t,x) \leq a_3(\gamma(t)|x|) \tag{2.3}$$

where $b \in C^0(\Re^+ \times (\Re^n \setminus \{0\}); \Re^+)$ is the function involved in (2.2).

**Remark 2.2:** If the ORCLF $V \in C^1(\Re^+ \times \Re^n; \Re^+)$ is $T$-periodic (namely, $V(t+T,x) = V(t,x)$ for certain $T > 0$ and for all $(t,x) \in \Re^+ \times \Re^n$), then (2.1) implies $a_1(M|x|) \leq V(t,x)$ for all $(t,x) \in \Re^+ \times \Re^n$, where $M := \min_{t \in [0,T]} \mu(t) > 0$. Consequently, existence of a $T$-periodic OCLF implies existence of a $T$-periodic SRCLF. Moreover, the existence of a time-invariant ORCLF implies the existence of a time-invariant SRCLF.

**Remark 2.3:** The small-control property in Definition 2.1 constitutes a time-varying version of the small-control property for the autonomous case [2,9,20].

**Remark 2.4:** Time-Varying ORCLFs have to be considered even for autonomous systems. It should be noticed that, in general, it is possible for an *autonomous* system (1.1) to possess a *time-varying* ORCLF satisfying the small-control property, although a *time-independent* ORCLF does not exists. Indeed, consider the elementary linear system $\dot{x}_1 = x_1$, $\dot{x}_2 = u$, with $u \in U = \Re$ and output $Y = x_2$. Obviously, this system is not feedback stabilizable to zero $0 \in \Re^2$ and therefore, according to [7], a time-invariant SRCLF does not exist. Neither a time-independent ORCLF exists, according to Remark 2.2 above. On the other hand it can be easily verified that the function $V(t,x) := \frac{1}{2}\exp(-4t)x_1^2 + \frac{1}{2}x_2^2$ is an ORCLF, which in addition satisfies the small-control property.

We next present certain stability concepts used in the present work. Consider the system

$$\dot{x} = f(t,d,x) \tag{2.4a}$$

$$\begin{aligned} Y &= H(t,x) \\ x &\in \Re^n, d \in D, Y \in \Re^k \end{aligned} \tag{2.4b}$$

where the mappings $f : \Re^+ \times D \times \Re^n \to \Re^n$, $H : \Re^+ \times \Re^n \to \Re^k$ are continuous with $f(t,d,0,0) = 0$, $H(t,0) = 0$ for all $(t,d) \in \Re^+ \times D$ and $D \subset \Re^l$ is compact. We assume that for every $(t_0, x_0, d) \in \Re^+ \times \Re^n \times M_D$ there exists $h \in (0, +\infty]$ and a unique solution $x(\cdot) = x(\cdot, t_0, x_0; d) : [t_0, t_0 + h) \to \Re^n$ of (2.4a) with $x(t_0) = x_0$.

**Definition 2.5:** We say that (2.4) is **Robustly Forward Complete (RFC),** if for every $T \geq 0$, $r \geq 0$ it holds that:

$$\sup\{|x(t_0+h, t_0, x_0; d)| \; ; \; |x_0| \leq r, t_0 \in [0,T], h \in [0,T], d(\cdot) \in M_D\} < +\infty \tag{2.5}$$

Clearly, the notion of robust forward completeness implies the standard notion of forward completeness, which simply requires that for every initial condition the solution of the system exists for all times greater than the initial time, or equivalently, the solutions of the system do not present finite escape time. Conversely, an extension of Proposition 5.1 in [16] to the time-varying case shows that every forward complete system (2.4) whose dynamics are locally Lipschitz with respect to $(t,x)$, uniformly in $d \in D$, is RFC. All output stability notions used in the present work will assume RFC.

We next provide the notion of (non-uniform in time) Robust Global Asymptotic Output Stability (RGAOS) (see [12,13]), which is a generalization of the notion of Robust Output Stability (see [22,23,26]). Let us denote by $Y(\cdot) = H(\cdot, x(\cdot, t_0, x_0; d))$ the output of (2.4) corresponding to input $d \in M_D$ and initial condition $x(t_0) = x_0$.

**Definition 2.6:** *Consider system (2.4) and suppose that is RFC. We say that (2.4) is **(non-uniformly in time) Robustly Globally Asymptotically Output Stable (RGAOS)** if it satisfies the following properties:*

**P1(Output Stability)** *For every $\varepsilon > 0$, $T \geq 0$, it holds:*

$$\sup\{|Y(t)|\,;\, t \geq t_0\,,\, |x_0| \leq \varepsilon\,,\, t_0 \in [0,T]\,,\, d(\cdot) \in M_D \} < +\infty$$
***(Robust Lagrange Output Stability)***

*and there exists a $\delta := \delta(\varepsilon, T) > 0$ such that*

$$|x_0| \leq \delta\,,\, t_0 \in [0,T] \Rightarrow |Y(t)| \leq \varepsilon\,,\, \forall t \geq t_0\,,\, \forall d(\cdot) \in M_D$$
***(Robust Lyapunov Output Stability)***

**P2(Uniform Output Attractivity on compact sets of initial data)** *For every $\varepsilon > 0$, $T \geq 0$ and $R \geq 0$, there exists a $\tau := \tau(\varepsilon, T, R) \geq 0$ such that*

$$|x_0| \leq R\,,\, t_0 \in [0,T] \Rightarrow |Y(t)| \leq \varepsilon\,,\, \forall t \geq t_0 + \tau\,,\, \forall d(\cdot) \in M_D$$

The notion of Uniform Robust Global Asymptotic Output Stability was originally given in [22,23] and is a special case of non-uniform in time RGAOS.

**Definition 2.7:** *Consider system (2.4) and suppose that is RFC. We say that (2.4) is **Uniformly Robustly Globally Asymptotically Output Stable (URGAOS)**, if it satisfies the following properties:*

**P1(Uniform Output Stability)** *For every $\varepsilon > 0$, it holds that*

$$\sup\{|Y(t)|\,;\, t \geq t_0\,,\, |x_0| \leq \varepsilon\,,\, t_0 \geq 0\,,\, d(\cdot) \in M_D \} < +\infty$$
***(Uniform Robust Lagrange Output Stability)***

*and there exists a $\delta := \delta(\varepsilon) > 0$ such that*

$$|x_0| \leq \delta\,,\, t_0 \geq 0 \Rightarrow |Y(t)| \leq \varepsilon\,,\, \forall t \geq t_0\,,\, \forall d(\cdot) \in M_D$$
***(Uniform Robust Lyapunov Output Stability)***

**P2(Uniform Output Attractivity on compact sets of initial states)** *For every $\varepsilon > 0$ and $R \geq 0$, there exists a $\tau := \tau(\varepsilon, R) \geq 0$ such that*

$$|x_0| \leq R\,,\, t_0 \geq 0 \Rightarrow |Y(t)| \leq \varepsilon\,,\, \forall t \geq t_0 + \tau\,,\, \forall d(\cdot) \in M_D$$

Obviously, for the case $H(t, x) = x$ the notions of RGAOS, URGAOS coincide with the notions of ***non-uniform in time Robust Global Asymptotic Stability (RGAS)*** as given in [11] and ***Uniform Robust Global Asymptotic Stability (URGAS)*** as given in [16], respectively. Also note that, if there exists $a \in K_\infty$ with $|x| \leq a(|H(t,x)|)$ for all $(t, x) \in \Re^+ \times \Re^n$, then (U)RGAOS implies (U)RGAS.

We are now in a position to state our main results.

**Theorem 2.8:** *Consider system (1.1) under hypotheses (H1-3) and assume that (1.1) admits an ORCLF which satisfies (2.1), (2.2) and the small-control property (2.3). Moreover, suppose that $\beta(t) \equiv 1$, where $\beta \in K^+$ is the function involved in (2.1). Then there exists a continuous mapping $K : \Re^+ \times \Re^n \to U$ with $K(t,0) = 0$ for all $t \geq 0$, being continuously differentiable with respect to $x \in \Re^n$ on the set $\Re^+ \times (\Re^n \setminus \{0\})$, such that for all $(t_0, x_0, d) \in \Re^+ \times \Re^n \times M_D$ the solution $x(\cdot)$ of the closed-loop system (1.1) with $u = K(t,x)$:*

$$\dot{x} = f(t, d, x, K(t, x)) \tag{2.6}$$

*with initial condition $x(t_0) = x_0 \in \Re^n$, corresponding to input $d \in M_D$ is unique and system (2.6) is URGAOS.*

**Theorem 2.9:** *Consider system (1.1) under hypotheses (H1-3) and assume that (1.1) admits an ORCLF which satisfies (2.1), (2.2). Then there exists a continuous mapping $K : \Re^+ \times \Re^n \to U$, with $K(t,0) = 0$ for all $t \geq 0$, which is continuously differentiable with respect to $x \in \Re^n$ on $\Re^+ \times \Re^n$, such that the closed-loop system (2.6) is RGAOS.*

It should be emphasized that the small-control property is not required for the validity of the result of Theorem 2.9. On the other hand, Theorem 2.9 cannot in general guarantee uniformity of solutions of the resulting closed-loop system (2.6) with respect to the initial time. Another advantage of Theorem 2.9 above is that the proposed feedback $K(t, x)$ is locally Lipschitz with respect to $x \in \Re^n$. The latter in conjunction with the converse Lyapunov theorem in [12] leads to the following result:

**Corollary 2.10:** *Consider system (1.1) under hypotheses (H1-3). The following statements are equivalent:*

(i) *System (1.1) admits an ORCLF which satisfies (2.1), (2.2).*

(ii) *There exists a continuous mapping $K : \Re^+ \times \Re^n \to U$, with $K(t,0) = 0$ for all $t \geq 0$, which is continuously differentiable with respect to $x \in \Re^n$ on $\Re^+ \times \Re^n$, such that the closed-loop system (1.1) with $u = K(t, x)$ is RGAOS.*

(iii) *System (1.1) admits an ORCLF which satisfies (2.1), (2.2) and the small-control property (2.3).*

The following example illustrates the nature of Theorem 2.9.

**Example 2.11:** Consider the following system

$$\begin{aligned} \dot{x} &= f(t, d, x) + \sum_{j=1}^{N} g_j(t, x) u^{k_j} \\ Y &= H(t, x) \\ x &\in \Re^n, d \in D, u \in \Re, Y \in \Re^k \end{aligned} \tag{2.7}$$

where $D \subset \Re^l$ is a compact set, $f : \Re^+ \times D \times \Re^n \to \Re^n$, $g_j : \Re^+ \times \Re^n \to \Re^n$ ($j = 0, ..., N$) are locally Lipschitz mappings with $f(t, d, 0) = 0$ for all $(t, d) \in \Re^+ \times D$ and $H : \Re^+ \times \Re^n \to \Re^k$ is a continuous mapping with $H(t, 0) = 0$ for all $t \geq 0$. Assume that

$$k_j, \; j = 1, ..., N \text{ are odd positive integers} \tag{2.8}$$

and there exist functions $V \in C^1(\Re^+ \times \Re^n; \Re^+)$, $a_1, a_2 \in K_\infty$, $\mu, \beta \in K^+$, $\rho \in C^0(\Re^+; \Re^+)$ being locally Lipschitz and positive definite, such that

$$a_1(|H(t, x)| + \mu(t)|x|) \leq V(t, x) \leq a_2(\beta(t)|x|), \; \forall (t, x) \in \Re^+ \times \Re^n \tag{2.9}$$

and in such a way that the following implication holds:

$$\sum_{j=1}^{N}\left(\frac{\partial V}{\partial x}(t,x)g_j(t,x)\right)^2 = 0 \Rightarrow \max_{d\in D}\left(\frac{\partial V}{\partial t}(t,x)+\frac{\partial V}{\partial x}(t,x)f(t,d,x)\right) \leq -2\rho(V(t,x)) \qquad (2.10)$$

We claim that $V \in C^1(\Re^+ \times \Re^n; \Re^+)$ is an ORCLF for system (2.7). Indeed, by exploiting (2.8) and implication (2.10) it follows that for every $(t,x) \in \Re^+ \times \Re^n$ there exists $u \in \Re$ such that

$$\max_{d\in D}\left(\frac{\partial V}{\partial t}(t,x)+\frac{\partial V}{\partial x}(t,x)f(t,d,x)+\sum_{j=1}^{N}\frac{\partial V}{\partial x}(t,x)g_j(t,x)u^{k_j}\right) \leq -2\rho(V(t,x)) \qquad (2.11)$$

From (2.11), compactness of $D \subset \Re^l$ and continuity of $f$, $g_j$ ($j=0,...,N$), it follows by applying standard partition of unity arguments, that there exists a function $b \in C^0(\Re^+ \times (\Re^n \setminus \{0\}); \Re^+)$ such that

$$\min_{\substack{|u|\leq b(t,x) \\ u\in\Re}} \max_{d\in D}\left(\frac{\partial V}{\partial t}(t,x)+\frac{\partial V}{\partial x}(t,x)f(t,d,x)+\sum_{j=1}^{N}\frac{\partial V}{\partial x}(t,x)g_j(t,x)u^{k_j}\right) \leq -\rho(V(t,x)) \qquad (2.12)$$

Hence, by (2.9) and (2.12) we may conclude that $V \in C^1(\Re^+ \times \Re^n; \Re^+)$ is an ORCLF for system (2.7). Consequently, according to statement of Theorem 2.9, there exists a continuous mapping $K: \Re^+ \times \Re^n \to U$, with $K(t,0)=0$ for all $t \geq 0$, which is continuously differentiable with respect to $x \in \Re^n$ on $\Re^+ \times \Re^n$, such that the closed-loop system (2.7) with $u = K(t,x)$ is RGAOS. ◁

## 3. Proofs of the Main Results

The proof of the main results of the present work is based on three lemmas below. Particularly, Lemma 3.1 is a preparatory result for the construction of the desired feedback stabilizer. It constitutes a time-varying extension of Lemma 2.7 in [7], but its constructive proof differs from the corresponding proof of the previously mentioned result.

**Lemma 3.1:** *Consider system (1.1) under hypotheses (H1-3) and assume that (1.1) admits an ORCLF which satisfies (2.1), (2.2). Then there exists a $C^1$ function $k:[0,1]\times\Re^+ \times (\Re^n \setminus \{0\}) \to U$ with*

$$k(0,t,x) = k(1,t,x) = 0 \qquad (3.1a)$$

$$\frac{\partial k}{\partial s}(0,t,x) = \frac{\partial k}{\partial t}(0,t,x) = 0 \; ; \; \frac{\partial k}{\partial x}(0,t,x) = 0 \qquad (3.1b)$$

$$\frac{\partial k}{\partial s}(1,t,x) = \frac{\partial k}{\partial t}(1,t,x) = 0 \; ; \; \frac{\partial k}{\partial x}(1,t,x) = 0 \qquad (3.1c)$$

*for all $t \geq 0$, $x \in \Re^n \setminus \{0\}$, and in such a way that:*

$$\frac{\partial V}{\partial t}(t,x) + \frac{\partial V}{\partial x}(t,x)\int_0^1 f(t,d(s),x,k(s,t,x))ds \leq -\frac{1}{2}\rho(V(t,x)) \qquad (3.2)$$

*for all $(t,x) \in \Re^+ \times (\Re^n \setminus \{0\})$, $d \in M_D$. Moreover, the following inequality holds for all $(t,x) \in \Re^+ \times (\Re^n \setminus \{0\})$:*

$$\max_{s\in[0,1]}|k(s,t,x)| \leq \tilde{b}(t,x) \qquad (3.3)$$

*where*

$$\tilde{b}(t,x) := \max\left\{ b(\tau, y) : \frac{2}{3}|x| \leq |y| \leq 2|x|, 0 \leq \tau \leq t+1 \right\} \tag{3.4}$$

**Proof of Lemma 3.1:** Let $\tilde{b} : \Re^+ \times \Re^n \to \Re^+$ as given by (3.4) that obviously is of class $C^0(\Re^+ \times \Re^n ; \Re^+)$ and let $\varphi : \Re^+ \times \Re^n \to [1,+\infty)$ be any smooth ($C^\infty$) function satisfying

$$\max_{\substack{|u| \leq \tilde{b}(t,x) \\ u \in U}} \max_{d \in D} \left( \frac{\partial V}{\partial t}(t,x) + \frac{\partial V}{\partial x}(t,x) f(t,d,x,u) \right) \leq \varphi(t,x), \ \forall (t,x) \in \Re^+ \times (\Re^n \setminus \{0\}) \tag{3.5}$$

Moreover, let $\varepsilon : \Re^+ \times (\Re^n \setminus \{0\}) \to (0,1)$ be a smooth function such that

$$0 < \varepsilon(t,x) \leq \frac{\rho(V(t,x))}{4(\rho(V(t,x) + \varphi(t,x)))}, \ \forall (t,x) \in \Re^+ \times (\Re^n \setminus \{0\}) \tag{3.6}$$

and define

$$\Psi(t,x,u) := \max_{d \in D} \left( \frac{\partial V}{\partial t}(t,x) + \frac{\partial V}{\partial x}(t,x) f(t,d,x,u) + \frac{3}{4}\rho(V(t,x)) \right), \ (t,x,u) \in \Re^+ \times \Re^n \times U \tag{3.7a}$$

$$\Psi(t,x,u) := \Psi(0,x,u), \ (t,x,u) \in (-1,0) \times \Re^n \times U \tag{3.7b}$$

By virtue of (2.2), continuity of $\Psi$ and compactness of $D \subset \Re^l$, it follows that for each $(t,x) \in (-1,+\infty) \times (\Re^n \setminus \{0\})$ there exist $u = u(t,x) \in U$ with $|u| \leq b(t,x)$ and $\delta = \delta(t,x) \in (0,1]$ with $\delta(t,x) \leq \frac{|x|}{2}$ such that

$$\Psi(\tau, y, u(t,x)) \leq 0, \ \forall (\tau, y) \in \left\{ (\tau, y) \in (-1,+\infty) \times \Re^n : |\tau - t| + |y - x| < \delta \right\} \tag{3.8}$$

Using (3.8) and standard partition of unity arguments, we can determine sequences $\{(t_i, x_i) \in (-1,+\infty) \times (\Re^n \setminus \{0\})\}_{i=1}^\infty$, $\{u_i \in U\}_{i=1}^\infty$, $\{\delta_i \in (0,1)\}_{i=1}^\infty$ with $|u_i| \leq b(\max(0,t_i), x_i)$ and $\delta_i = \delta(t_i, x_i) \leq \frac{|x_i|}{2}$ associated with a sequence of open sets $\{\Omega_i\}_{i=1}^\infty$ with

$$\Omega_i \subseteq \left\{ (\tau, y) \in (-1,+\infty) \times \Re^n : |\tau - t_i| + |y - x_i| < \delta_i \right\} \tag{3.9a}$$

forming a locally finite open covering of $(-1,+\infty) \times (\Re^n \setminus \{0\})$ and in such a way that:

$$\Psi(\tau, y, u_i) \leq 0, \ \forall (\tau, y) \in \Omega_i \tag{3.9b}$$

Also, a family of smooth functions $\{\theta_i\}_{i=1}^\infty$ with $\theta_i(t,x) \geq 0$ for all $(t,x) \in (-1,+\infty) \times (\Re^n \setminus \{0\})$ can be determined with

$$supp \ \theta_i \subseteq \Omega_i \tag{3.9c}$$

$$\sum_{i=1}^\infty \theta_i(t,x) = 1, \ \forall (t,x) \in (-1,+\infty) \times (\Re^n \setminus \{0\}) \tag{3.9d}$$

Next define recursively the following mappings for each $(t,x) \in \Re^+ \times (\Re^n \setminus \{0\})$:

$$T_i(t,x) = T_{i-1}(t,x) + \theta_i(t,x), \ i \geq 1; \ T_0(t,x) = 0; \ (t,x) \in \Re^+ \times (\Re^n \setminus \{0\}) \tag{3.10}$$

Notice that definition (3.10) implies $T_n(x) = \sum_{i=1}^{n} \theta_i(t,x)$ for all $n \geq 1$. Since the open sets $\{\Omega_i\}_{i=1}^{\infty}$ form a locally finite open covering of $\Re^+ \times (\Re^n \setminus \{0\})$, it follows from (3.9c) and (3.10) that for every $(t,x) \in \Re^+ \times (\Re^n \setminus \{0\})$ there exists $m = m(t,x) \in \{1,2,3,...\}$ such that

$$T_i(t,x) = 1 \text{ for } i \geq m \tag{3.11}$$

We define the index set

$$J(t,x) := \{ j \in \{1,2,3,...\} : \theta_j(t,x) > 0 \} \tag{3.12}$$

which by virtue of (3.11) is a non-empty finite set. It follows from definitions (3.10) and (3.12) that

$$\bigcup_{j \in J(t,x)} [T_{j-1}(t,x), T_j(t,x)) = [0,1), \ \forall (t,x) \in [0,1] \times \Re^+ \times (\Re^n \setminus \{0\}) \tag{3.13}$$

Let $h : \Re \to [0,1]$ be any smooth non-decreasing function with

$$h(s) = 0 \text{ for } s \leq 0 \text{ and } h(s) = 1 \text{ for } s \geq 1 \tag{3.14a}$$

and let

$$g_j(t,x) := \frac{1}{2}\theta_j(t,x) + \frac{1}{2}\left(\varepsilon(t,x)2^{-j-1} - \theta_j(t,x)\right) h\left(\frac{\theta_j(t,x) - \varepsilon(t,x)2^{-j-2}}{\varepsilon(t,x)2^{-j-2}}\right), \ j \in \{1,2,3,...\} \tag{3.14b}$$

where $\varepsilon(\cdot,\cdot)$ is the function defined by (3.6). Notice that according to (3.14a,b) it holds:

$$\min\left\{\varepsilon(t,x)2^{-j-2}, \frac{1}{2}\theta_j(t,x)\right\} \leq g_j(t,x) \leq \min\left\{\varepsilon(t,x)2^{-j-2}, \theta_j(t,x)\right\} \tag{3.15a}$$

$$g_j(t,x) = \varepsilon(t,x)2^{-j-2} \text{ for } \theta_j(t,x) \geq \varepsilon(t,x)2^{-j-1} \tag{3.15b}$$

We define the following map $[0,1] \times \Re^+ \times (\Re^n \setminus \{0\}) \ni (s,t,x) \to k(s,t,x) \in \Re^m$:

$$k(s,t,x) = u_j \left(\frac{g_j(t,x)}{\varepsilon(t,x)2^{-j-2}}\right)^2 h\left(\frac{s - T_{j-1}(t,x) - \frac{1}{5}g_j(t,x)}{\frac{1}{5}g_j(t,x)}\right), \text{ for } s \in \left[T_{j-1}(t,x), T_{j-1}(t,x) + \frac{1}{2}\theta_j(t,x)\right), \ j \in J(t,x)$$

$$\tag{3.16a}$$

$$k(s,t,x) = u_j \left(\frac{g_j(t,x)}{\varepsilon(t,x)2^{-j-2}}\right)^2 h\left(\frac{T_j(t,x) - \frac{1}{5}g_j(t,x) - s}{\frac{1}{5}g_j(t,x)}\right), \text{ for } s \in \left[T_j(t,x) - \frac{1}{2}\theta_j(t,x), T_j(t,x)\right), \ j \in J(t,x)$$

$$\tag{3.16b}$$

$$k(1,t,x) = 0 \tag{3.16c}$$

Notice that because of (3.13), $k(\cdot,\cdot,\cdot)$ is well defined for all $(s,t,x) \in [0,1] \times \Re^+ \times (\Re^n \setminus \{0\})$. Furthermore, according to definition (3.16), hypothesis (H2) guarantees that $k(\cdot,\cdot,\cdot)$ takes values in $U \subseteq \Re^m$ and is continuously differentiable on the region $\left( \bigcup_{j \in J(t,x)} (T_{j-1}(t,x), T_j(t,x)) \right) \times \Re^+ \times (\Re^n \setminus \{0\})$. Furthermore, it holds that

$$\frac{\partial k}{\partial s}(s,t,x) \to 0, \ \frac{\partial k}{\partial s}(s,t,x) \to 0, \ \frac{\partial k}{\partial x}(s,t,x) \to 0 \text{ as } s \to T_j(t,x) \text{ for all } j \in \{0,1,2,3,...\} \quad (3.17)$$

Next, we show that $k(\cdot,\cdot,\cdot)$ is continuously differentiable on the whole region $[0,1] \times \Re^+ \times (\Re^n \setminus \{0\})$ and simultaneously that (3.1b,c,d) are fulfilled. We distinguish the following cases:

Case 1: Let $s \in (0,1)$, $(t,x) \in \Re^+ \times (\Re^n \setminus \{0\})$ and suppose that there exists a positive integer $p$ with $s = T_p(t,x)$. Then, there exist positive integers $m, l$ with $l \leq p \leq m$ in such a way that

$$\theta_{m+1}(t,x) > 0, \ \theta_l(t,x) > 0 \quad (3.18a)$$

$$s = T_m(t,x) = ... = T_l(t,x) > 0 \quad (3.18b)$$

Equality (3.18b) in conjunction with definition (3.10) means

$$\theta_m(t,x) = ... = \theta_{l+1}(t,x) = 0, \text{ if } m \geq l+1 \quad (3.19)$$

Notice that definition (3.16a) and (3.18a) imply that in our case it holds

$$k(s,t,x) = 0 \quad (3.20)$$

By taking into account continuity of the mappings $g_l, g_{m+1}, T_l, T_m$ and (3.15a), it follows that there exists $\delta > 0$ such that

$$s' \in \left( T_l(\tau, y) - \frac{1}{5} g_l(\tau, y), T_m(\tau, y) + \frac{1}{5} g_{m+1}(\tau, y) \right)$$
$$\forall (s', \tau, y) \in [0,1] \times \Re^+ \times (\Re^n \setminus \{0\}) \text{ with } |s' - s| + |\tau - t| + |y - x| < \delta \quad (3.21)$$

By virtue of definition (3.16a,b), (3.20) and (3.21) it follows that for every $(s', \tau, y) \in [0,1] \times \Re^+ \times (\Re^n \setminus \{0\})$ with $|s' - s| + |\tau - t| + |y - x| < \delta$ it holds:

$$|k(s',\tau,y) - k(s,t,x)| \leq \max_{v=l+1,...,m} |u_v| \left( \frac{g_v(\tau,y)}{\varepsilon(\tau,y) 2^{-v-2}} \right)^2, \text{ if } m \geq l+1 \quad (3.22a)$$

$$k(s',\tau,y) = k(s,t,x), \text{ if } m = l \quad (3.22b)$$

If $m \geq l+1$, then by (3.15a) and (3.19) we also get $g_v(t,x) = 0$ for $v = l+1,...,m$, hence, since the mappings $g_v$ are continuously differentiable, there exists a constant $L > 0$ such that

$$\frac{\max_{v=l+1,...,m} g_v(\tau,y)}{\varepsilon(\tau,y)} \leq L|\tau - t| + L|y - x|, \ \forall (\tau,y) \in \Re^+ \times \Re^n \setminus \{0\} \text{ with } |\tau - t| + |y - x| < \delta \quad (3.23)$$

It turns out from (3.22a,b) and (3.23) that

$$|k(s',\tau,y) - k(s,t,x)| \leq L' \left( |\tau - t|^2 + |y - x|^2 \right) \quad (3.24)$$

for certain constant $L' > 0$ and for $|s' - s| + |\tau - t| + |y - x| < \delta$. We conclude from (3.24) that the derivatives of $k(\cdot, \cdot, \cdot)$ exist for $s = T_p(t, x)$ and it holds that $\frac{\partial k}{\partial s}(s, t, x) = \frac{\partial k}{\partial t}(s, t, x) = 0$ and $\frac{\partial k}{\partial x}(s, t, x) = 0$ for $s = T_p(t, x)$. The latter in conjunction with (3.17) implies that $k(\cdot, \cdot, \cdot)$ is continuously differentiable in a neighborhood of $(s, t, x)$ with $s = T_p(t, x) \in (0,1)$.

Case 2: Let $s = 0$, $(t, x) \in \Re^+ \times (\Re^n \setminus \{0\})$ and suppose that there exists an integer $p \geq 0$ with $s = T_p(t, x) = 0$. Clearly, there exists an integer $m \geq p$ such that

$$\theta_{m+1}(t, x) > 0 \tag{3.25a}$$

$$s = T_m(t, x) = ... = T_0(t, x) = 0 \tag{3.25b}$$

(note again that equality (3.25b) means that $\theta_m(t, x) = ... = \theta_1(t, x) = 0$ for the case $m > 0$). By virtue of definition (3.16a) it holds that $k(s, t, x) = 0$ and continuity of the mappings $T_m$ and $g_{m+1}$ implies that there exists $\delta > 0$ such that $s' \in \left[0, T_m(\tau, y) + \frac{1}{5} g_{m+1}(\tau, y)\right)$ for all $(s', \tau, y) \in [0,1] \times \Re^+ \times (\Re^n \setminus \{0\})$ with $|s' - s| + |\tau - t| + |y - x| < \delta$. Then as in Case 1, it follows by (3.16) that

$$|k(s', \tau, y) - k(s, t, x)| \leq \max_{v=1,...,m} |u_v| \left(\frac{g_v(\tau, y)}{\varepsilon(\tau, y) 2^{-v-2}}\right)^2, \text{ if } m > 0 \tag{3.26a}$$

$$k(s', \tau, y) = k(s, t, x), \text{ if } m = 0 \tag{3.26b}$$

for every $|s' - s| + |\tau - t| + |y - x| < \delta$, from which we get the desired conclusion, namely, that that $k(\cdot, \cdot, \cdot)$ is continuously differentiable in a neighborhood of $(0, t, x)$ and further (3.1b) holds.

Case 3: Let $s = 1$, $(t, x) \in \Re^+ \times (\Re^n \setminus \{0\})$ and let $p$ a positive integer with $s = T_p(t, x)$. Let $\{\Omega_i\}_{i=1}^\infty$ be the locally finite open covering of $(-1, +\infty) \times (\Re^n \setminus \{0\})$ and the associated sequence of functions $\{\theta_i\}_{i=1}^\infty$ in such a way that (3.9a,b,c,d) hold. Let $N \subset (-1, +\infty) \times (\Re^n \setminus \{0\})$ be a neighborhood containing $(t, x)$ which intersects only a finite number of the open sets $\{\Omega_i\}_{i=1}^\infty$ (see [10]). Consequently, by (3.9d) there exists an integer $m > 1$ such that $N \cap \Omega_i = \emptyset$ for all $i > m$ and $\theta_i(\tau, y) = 0$ for all $i > m$, $(\tau, y) \in N$. Clearly, there exists $l \in \{1,...,m\}$ with

$$\theta_l(t, x) > 0 \tag{3.27a}$$

$$s = T_l(t, x) = ... = T_m(t, x) = 1 \tag{3.27b}$$

Without loss of generality we may assume that $m > l$. By virtue of definition (3.16c) we have $k(1, t, x) = 0$ and continuity of the mappings $T_l$, $g_l$, asserts existence of a constant $\delta > 0$ such that

$$s' \in \left(T_l(\tau, y) - \frac{1}{5} g_l(\tau, y), 1\right] \text{ and } (\tau, y) \in N ;$$
$$\forall (s', \tau, y) \in [0,1] \times \Re^+ \times (\Re^n \setminus \{0\}) \text{ with } |s' - 1| + |\tau - t| + |y - x| < \delta \tag{3.28}$$

Using (3.16) and (3.28) we get

$$|k(s', \tau, y) - k(1, t, x)| \leq \max_{v=l+1,...,m} |u_v| \left(\frac{g_v(\tau, y)}{\varepsilon(\tau, y) 2^{-v-2}}\right)^2$$

from which it follows that (3.24) holds for all $|s'-1|+|\tau-t|+|y-x|<\delta$ and for certain constant $L'>0$. This implies that the derivatives of $k(\cdot,\cdot,\cdot)$ exist for $s=1$ and particularly, (3.1c) holds. The latter in conjunction with (3.17) implies that $k(\cdot,\cdot,\cdot)$ is continuously differentiable in a neighborhood of $(1,t,x)$.

We next establish (3.3). By virtue of (3.14a), (3.15a) and definition (3.16) we have $\max_{s\in[0,1]}|k(s,t,x)|\leq \max_{j\in J(t,x)}|u_j|$, $J(t,x)$ being the index set defined by (3.12). For every $j\in J(t,x)$ there exist $(t_j,x_j)\in(-1,+\infty)\times(\Re^n\setminus\{0\})$ with

$$|u_j|\leq b(\max(0,t_j),x_j) \tag{3.29}$$

for which $(t,x)\in\Omega_j$ and in such a way that (3.9a) holds with $i=j$. Since $\delta_j=\delta(t_j,x_j)\leq\min\left\{1,\frac{|x_j|}{2}\right\}$, it follows that, when $|t-t_j|+|x-x_j|<\delta_j$, it holds that $\frac{2}{3}|x|\leq|x_j|\leq 2|x|$ and $t-1\leq t_j\leq t+1$. The latter in conjunction with (3.25) and definition (3.4) of $\tilde{b}(\cdot,\cdot)$ implies (3.3). Finally, we establish (3.2). Notice that by (3.12), (3.14a), (3.15b), (3.16a,b), for any $(t,x)\in\Re^+\times\Re^n\setminus\{0\}$ and integer $j\in J(t,x)$ it holds:

$$k(s,t,x)=u_j,\ \forall s\in\left[T_{j-1}(t,x)+\frac{2}{5}g_j(t,x), T_j(t,x)-\frac{2}{5}g_j(t,x)\right] \text{ when } \theta_j(t,x)\geq\varepsilon(t,x)2^{-j-1} \tag{3.30}$$

hence, the set $I_{(t,x)}:=\{s\in[0,1]:k(s,t,x)\neq u_j,\ j\in J(t,x)\}$ has Lebesgue measure, say $|I_{(t,x)}|$, satisfying:

$$|I_{(t,x)}|\leq \sum_{j\in J(t,x)}\varepsilon(t,x)2^{-j-1}\leq\varepsilon(t,x) \tag{3.31}$$

Then for any $d\in M_D$ it follows by virtue of (3.7a), (3.9b) and (3.31) that

$$\frac{\partial V}{\partial t}(t,x)+\frac{\partial V}{\partial x}(t,x)\int_0^1 f(t,d(s),x,k(s,t,x))ds\leq$$
$$\leq -\frac{3}{4}(1-\varepsilon(t,x))\rho(V(t,x))+\varepsilon(t,x)\max_{\substack{|u|\leq\tilde{b}(t,x)\\u\in U}}\max_{d\in D}\frac{\partial V}{\partial t}(t,x)+\frac{\partial V}{\partial x}(t,x)f(t,d,x,u)$$

Inequalities (3.5), (3.6) in conjunction with the above inequality imply (3.2) and the proof is complete. ◁

The next lemmas (Lemma 3.2 and 3.3) constitute key results of the rest analysis and generalize Lemmas 2.8, 2.9 in [7]. Their proofs are based on certain appropriate generalizations of the technique employed in [7].

**Lemma 3.2:** *Consider system (1.1) under the same hypotheses with those imposed in Lemma 3.1. For every pair of sets $r=\{r_i:i\in Z\}$, $a=\{a_i:i\in Z\}$ with $r_i>0$, $a_i>0$,*

$$r_i+2a_i<r_{i+1}-2a_{i+1}\text{ for all }i\in Z \tag{3.32}$$

$$\lim_{i\to+\infty}r_i=+\infty,\ \lim_{i\to-\infty}r_i=0 \tag{3.33}$$

*there exists a continuous mapping $k_{r,a}:\Re^+\times(\Re^n\setminus\{0\})\to U$, being continuously differentiable with respect to $x\in\Re^n\setminus\{0\}$ with*

$$k_{r,a}(j,x) = 0 \text{ and } \frac{\partial k_{r,a}}{\partial x}(j,x) = 0 \text{ for all } (x,j) \in \left(\Re^n \setminus \{0\}\right) \times Z^+ \quad (3.34)$$

$$\left|k_{r,a}(t,x)\right| \leq \tilde{b}(t,x), \ \forall (t,x) \in \Re^+ \times (\Re^n \setminus \{0\}) \quad (3.35)$$

where $\tilde{b}(\cdot,\cdot)$ is defined by (3.4), and in such a way that the following property holds for all $(t_0, x_0, d, i) \in \Re^+ \times \left(\Re^n \setminus \{0\}\right) \times M_D \times Z$:

$$V(t_0, x_0) \in [r_{i-1}, r_i] \implies V(t, x(t, t_0, x_0; d)) \leq r_i + \frac{5}{2} a_i, \text{ for all } t \in [t_0, \min([t_0]+1, t_{\max})) \quad (3.36)$$

where $x(\cdot, t_0, x_0; d)$ denotes the unique solution of

$$\dot{x} = f(t, d, x, k_{r,a}(t,x)), (t,x) \in \Re^+ \times (\Re^n \setminus \{0\}) \quad (3.37)$$

with initial condition $x(t_0) = x_0 \in \Re^n \setminus \{0\}$, corresponding to $d \in M_D$, $t_{\max} := t_{\max}(t_0, x_0, d) > t_0$ denotes its maximal existence time. Moreover, for each $(x_0, j, i) \in \left(\Re^n \setminus \{0\}\right) \times Z^+ \times Z$, there exists a positive integer $N \geq 2$ such that

$$V(j, x_0) \leq r_i - 2a_i \implies V\left(j + \frac{s}{N}, x\left(j + \frac{s}{N}, j, x_0; d\right)\right) \leq \max\left(r_{i-1} + 2a_{i-1}, V(j, x_0) - \frac{s}{N} \mu_i\right),$$

$$\text{for all } d \in M_D, \ s \in \{0, 1, \ldots, N\} \text{ with } j + \frac{s}{N} < t_{\max} \quad (3.38)$$

where

$$\mu_i := \frac{1}{4} \min\{\rho(s) : s \in [r_{i-1}, r_i]\} \quad (3.39)$$

**Proof of Lemma 3.2:** Let $k : [0,1] \times \Re^+ \times \left(\Re^n \setminus \{0\}\right) \to U$ be a $C^1$ function which satisfies (3.1), (3.2), (3.3) and whose existence is guaranteed by Lemma 3.1. Let $i \in Z$, $j \in Z^+$, define

$$\Omega_{i,j} = \{(t,x) \in [j, j+1] \times \Re^n : V(t,x) \in [r_{i-1}, r_i]\} \quad (3.40)$$

$$\rho_i := \min(a_{i-2}, a_{i-1}, a_i, a_{i+1}) \quad (3.41)$$

and let $\delta_{i,j} > 0$ satisfying:

$$\left|\frac{\partial V}{\partial t}(t,x) - \frac{\partial V}{\partial t}(t_0, x_0)\right| + \left|\frac{\partial V}{\partial x}(t,x) - \frac{\partial V}{\partial x}(t_0, x_0)\right| \max\left\{\left|f(t,x,d,u)\right| : d \in D, u \in U, |u| \leq \tilde{b}(t,x)\right\}$$

$$+ \left|\frac{\partial V}{\partial x}(t_0, x_0)\right| \max\left\{\left|f(t, d, x, k(s, t_0, x)) - f(t_0, d, x_0, k(s, t_0, x_0))\right| : s \in [0,1], d \in D\right\} \leq \frac{1}{4} \rho(V(t_0, x_0))$$

$$\forall (t_0, x_0) \in \Omega_{i,j}, \ \forall (t,x) \in \Re^+ \times (\Re^n \setminus \{0\}) \text{ with } t \in [t_0, t_0 + \delta_{i,j}], |x - x_0| \leq \delta_{i,j} \quad (3.42)$$

Also, let $N_{i,j} \in Z^+$ with $N_{i,j} \geq 2$ be a family of integers which satisfies the following inequalities:

$$4 \max\left\{\left|\frac{\partial V}{\partial t}(t,x) + \frac{\partial V}{\partial x}(t,x) f(t,d,x,u)\right| : t \in [j, j+2], V(t,x) \in [r_{i-3}, r_{i+2}], d \in D, u \in U, |u| \leq \tilde{b}(t,x)\right\} \leq \rho_i N_{i,j}$$

$$(3.43)$$

$$2 + 2\max\left\{|f(t,d,x,u)| : t \in [j, j+2], V(t,x) \in [r_{i-3}, r_{i+2}], d \in D, u \in U, |u| \le \tilde{b}(t,x)\right\} \le \delta_{i,j} N_{i,j} \tag{3.44}$$

Consider next a smooth non-decreasing function $h : \Re \to [0,1]$ with $h(s) = 0$ for $s \le 0$ and $h(s) = 1$ for $s \ge 1$ and define the desired $k_{r,a} : \Re^+ \times (\Re^n \setminus \{0\}) \to U$ as follows:

$$k_{r,a}(t,x) := \begin{cases} h\left(2\dfrac{V(t,x) - r_{i-1}}{\min(a_{i-1}, a_i)}\right) k\left(N_{i,j}(t-j) - l, j + \dfrac{l}{N_{i,j}}, x\right), & V(t,x) \in \left[r_{i-1}, \dfrac{r_i + r_{i-1}}{2}\right) \\ h\left(2\dfrac{r_i - V(t,x)}{\min(a_{i-1}, a_i)}\right) k\left(N_{i,j}(t-j) - l, j + \dfrac{l}{N_{i,j}}, x\right), & V(t,x) \in \left[\dfrac{r_i + r_{i-1}}{2}, r_i\right) \end{cases}$$

$$(t,x) \in \Omega_{i,j}, \ t \in \left[j + \dfrac{l}{N_{i,j}}, j + \dfrac{l+1}{N_{i,j}}\right) \text{ for some } l \in \{0,1,\ldots, N_{i,j} - 1\} \tag{3.45}$$

Obviously, (3.35) is a consequence of (3.3), (3.4) and (3.45). Moreover, by taking into account (3.1), (3.32), it follows that $k_{r,a}(\cdot, \cdot)$ above is continuous, continuously differentiable with respect to $x \in \Re^n \setminus \{0\}$ and satisfies

$$k_{r,a}(j,x) = 0, \ \frac{\partial k_{r,a}}{\partial x}(j,x) = 0, \ \forall (x,j) \in (\Re^n \setminus \{0\}) \times Z^+ \tag{3.46}$$

Let $(x_0, d) \in (\Re^n \setminus \{0\}) \times M_D$ and $t_0 \in \left[j + \dfrac{l}{N_{i,j}}, j + \dfrac{l+1}{N_{i,j}}\right)$ for some $l \in \{0, 1, \ldots, N_{i,j} - 1\}$ with

$$V(t_0, x_0) \in [r_{i-2}, r_{i+1}] \tag{3.47}$$

Then by (3.43), (3.44) and (3.47) it can be easily established that for all $t \in \left[t_0, j + \dfrac{l+1}{N_{i,j}}\right]$ it holds:

$$t - t_0 + |x(t, t_0, x_0; d) - x_0| \le \delta_{i,j} \tag{3.48a}$$

$$V(t_0, x_0) - \frac{1}{2}\min(a_{i-1}, a_i) \le V(t, x(t, t_0, x_0; d)) \le V(t_0, x_0) + \frac{1}{2}\min(a_{i-1}, a_i) \tag{3.48b}$$

Indeed, suppose on the contrary that there exist $(x_0, d) \in (\Re^n \setminus \{0\}) \times M_D$, $t_0 \in \left[j + \dfrac{l}{N_{i,j}}, j + \dfrac{l+1}{N_{i,j}}\right)$ for some $l \in \{0, 1, \ldots, N_{i,j} - 1\}$ satisfying (3.47) and $\bar{t} \in \left[t_0, j + \dfrac{l+1}{N_{i,j}}\right]$ such that either (3.48a) or (3.48b) does not hold and consider the closed set

$$A := \left\{\tau \in \left[t_0, j + \frac{l+1}{N_{i,j}}\right] : \max\left\{\frac{2|V(\tau, x(\tau, t_0, x_0; d)) - V(t_0, x_0)|}{\min(a_{i-2}, a_{i-1}, a_i, a_{i+1})}, \frac{\tau - t_0 + |x(\tau, t_0, x_0; d) - x_0|}{\delta_{i,j}}\right\} \ge 1\right\}$$

Notice that, since $\bar{t} \in A$, the set $A$ is non-empty. Let $t_1 := \min A$. Clearly, since $t_0 \notin A$, it holds that $t_1 > t_0$. Definition of the set $A$ above, (3.32) and (3.47) imply that $V(\tau, x(\tau, t_0, x_0; d)) \in [r_{i-3}, r_{i+2}]$ for every $\tau \in [t_0, t_1)$. It follows from (3.35), (3.43), (3.44) that

$$\left|\frac{d}{d\tau}V(\tau, x(\tau, t_0, x_0; d))\right| \le \frac{1}{4}\rho_i N_{i,j} \text{ and } 2 + 2|\dot{x}(\tau)| \le \delta_{i,j} N_{i,j}, \text{ a.e. for } \tau \in [t_0, t_1)$$

which in conjunction with definition (3.41) and the fact that $\tau - t_0 \leq \frac{1}{N_{i,j}}$ imply that for all $\tau \in [t_0, t_1]$ we would have:

$$|V(\tau, x(\tau, t_0, x_0; d)) - V(t_0, x_0)| \leq \int_{t_0}^{\tau} \left| \frac{d}{ds} V(s, x(s, t_0, x_0; d)) \right| ds \leq \frac{1}{4} \min(a_{i-2}, a_{i-1}, a_i, a_{i+1}) ;$$

$$\tau - t_0 + |x(\tau, t_0, x_0; d) - x_0| \leq \tau - t_0 + \int_{t_0}^{\tau} |\dot{x}(s)| ds \leq \frac{1}{2} \delta_{i,j}$$

The previous inequalities for $\tau = t_1$ are in contradiction with the fact that $t_1 \in A$.

In order to establish properties (3.36) and (3.38), we first need the following properties:

**Property P1:** *Let $d \in M_D$ and let*

$$t_0 = j + \frac{l}{N_{i,j}} \quad , l \in \{0, 1, ..., N_{i,j} - 1\} \tag{3.49a}$$

$$V(t_0, x_0) \in [r_{i-1} + a_{i-1}, r_i - 2a_i] \tag{3.49b}$$

*Then the following inequality is fulfilled:*

$$0 < V\left(j + \frac{l+1}{N_{i,j}}, x\left(j + \frac{l+1}{N_{i,j}}, j + \frac{l}{N_{i,j}}, x_0; d\right)\right) \leq V\left(j + \frac{l}{N_{i,j}}, x_0\right) - \frac{1}{4N_{i,j}} \rho\left(V\left(j + \frac{l}{N_{i,j}}, x_0\right)\right) \tag{3.50}$$

*Proof of P1:* Using (3.48b) and definition (3.45) it follows:

$$k_{r,a}(t, x(t, t_0, x_0; d)) = k\left(N_{i,j}(t-j) - l, j + \frac{l}{N_{i,j}}, x(t, t_0, x_0; d)\right), \forall t \in \left[j + \frac{l}{N_{i,j}}, j + \frac{l+1}{N_{i,j}}\right] \tag{3.51}$$

For convenience let us denote here $h := \frac{1}{N_{i,j}}$, $x(\cdot) = x(\cdot, t_0, x_0; d)$ and $\tilde{d}(t) := d(t_0 + ht)$ (notice that $\tilde{d} \in M_D$). From (3.51) we have:

$$V(t_0 + h, x(t_0 + h)) - V(t_0, x_0) = \int_{t_0}^{t_0+h} \left[ \frac{\partial V}{\partial t}(\tau, x(\tau)) + \frac{\partial V}{\partial x}(\tau, x(\tau)) f\left(\tau, d(\tau), x(\tau), k\left(\frac{\tau - t_0}{h}, t_0, x(\tau)\right)\right) \right] d\tau$$

$$= h \int_0^1 \left[ \frac{\partial V}{\partial t}(t_0 + hs, x(t_0 + hs)) + \frac{\partial V}{\partial x}(t_0 + hs, x(t_0 + hs)) f(t_0 + hs, d(t_0 + hs), x(t_0 + hs), k(s, t_0, x(t_0 + hs))) \right] ds$$

$$= h \int_0^1 \left[ \frac{\partial V}{\partial t}(t_0, x_0) + \frac{\partial V}{\partial x}(t_0, x_0) f(t_0, \tilde{d}(s), x_0, k(s, t_0, x_0)) \right] ds +$$

$$+ h \int_0^1 \left[ \frac{\partial V}{\partial t}(t_0 + hs, x(t_0 + hs)) - \frac{\partial V}{\partial t}(t_0, x_0) \right] ds +$$

$$+ h \int_0^1 \left[ \frac{\partial V}{\partial x}(t_0 + hs, x(t_0 + hs)) - \frac{\partial V}{\partial x}(t_0, x_0) \right] f(t_0 + hs, \tilde{d}(s), x(t_0 + hs), k(s, t_0, x(t_0 + hs))) ds +$$

$$+ h \int_0^1 \frac{\partial V}{\partial x}(t_0, x_0) \left[ f(t_0 + hs, \tilde{d}(s), x(t_0 + hs), k(s, t_0, x(t_0 + hs))) - f(t_0, \tilde{d}(s), x_0, k(s, t_0, x_0)) \right] ds$$

$$\tag{3.52}$$

Using (3.2), (3.3), (3.42), (3.48), (3.49) and (3.52) we get the desired (3.50) and the proof of P1 is complete.

The next property is a consequence of P1:

**Property P2**: *Suppose that*

$$0 < V\left(j + \frac{l}{N_{i,j}}, x_0\right) \leq r_i - 2a_i \text{ for some } l \in \{0,1,\ldots, N_{i,j} - 1\} \tag{3.53}$$

*and assume that the solution of (3.37) with initial condition* $x\left(j + \frac{l}{N_{i,j}}\right) = x_0 \in \Re^n \setminus \{0\}$, *corresponding to some* $d \in M_D$ *exists for* $t \in \left[j + \frac{l}{N_{i,j}}, j + \frac{l+1}{N_{i,j}}\right]$. *Then*

$$0 < V\left(j + \frac{l+1}{N_{i,j}}, x\left(j + \frac{l+1}{N_{i,j}}, j + \frac{l}{N_{i,j}}, x_0; d\right)\right) \leq \max\left\{r_{i-1} + 2a_{i-1}, V\left(j + \frac{l}{N_{i,j}}, x_0\right) - \frac{1}{N_{i,j}}\mu_i\right\} \tag{3.54}$$

*where* $\mu_i > 0$ *is defined by (3.39).*

*Proof of P2:* Obviously, the desired (3.54) is a consequence of (3.50), provided that (3.49b) is fulfilled. Consider the remaining case

$$0 < V\left(j + \frac{l}{N_{i,j}}, x_0\right) \leq r_{i-1} + a_{i-1} \tag{3.55}$$

We show by contradiction that, when (3.55) holds, then $0 < V\left(j + \frac{l+1}{N_{i,j}}, x\left(j + \frac{l+1}{N_{i,j}}, j + \frac{l}{N_{i,j}}, x_0; d\right)\right) \leq r_{i-1} + 2a_{i-1}$.

Indeed, suppose on the contrary that

$$V\left(j + \frac{l+1}{N_{i,j}}, x\left(j + \frac{l+1}{N_{i,j}}, j + \frac{l}{N_{i,j}}, x_0; d\right)\right) > r_{i-1} + 2a_{i-1} \tag{3.56}$$

Then, there would exist $t_1 \in \left(j + \frac{l}{N_{i,j}}, j + \frac{l+1}{N_{i,j}}\right)$ in such a way that $V\left(t_1, x\left(t_1, j + \frac{l}{N_{i,j}}, x_0; d\right)\right) = r_{i-1} + \frac{3a_{i-1}}{2}$.

Using (3.48b) the latter implies $0 < V\left(j + \frac{l+1}{N_{i,j}}, \phi\left(j + \frac{l+1}{N_{i,j}}, j + \frac{l}{N_{i,j}}, x_0; d\right)\right) \leq r_{i-1} + 2a_{i-1}$ which contradicts (3.56), and the proof of the P2 is complete.

The following property is a direct consequence of property P2 and (3.32):

**Property P3:** *Suppose that (3.53) holds and assume again that the solution of (3.37) with initial condition* $x\left(j + \frac{l}{N_{i,j}}\right) = x_0 \in \Re^n \setminus \{0\}$ *corresponding to some* $d \in M_D$ *exists for* $t \in \left[j + \frac{l}{N_{i,j}}, j + \frac{l+s}{N_{i,j}}\right]$, *for certain* $s \in \{0,1,2,\ldots, N_{i,j} - l\}$. *Then*

$$0 < V\left(j + \frac{l+s}{N_{i,j}}, x\left(j + \frac{l+s}{N_{i,j}}, j + \frac{l}{N_{i,j}}, x_0; d\right)\right) \leq \max\left\{r_{i-1} + 2a_{i-1}, V\left(j + \frac{l}{N_{i,j}}, x_0\right) - \frac{s}{N_{i,j}}\mu_i\right\} \tag{3.57}$$

The desired (3.38) follows from property P3 with $l = 0$, $N = N_{i,j}$. We next proceed with the proof of (3.36). Combining property P3 with (3.48b) we obtain:

**Property P4:** *If (3.53) is fulfilled then*

$$0 < V\left(t, x\left(t, j + \frac{l}{N_{i,j}}, x_0; d\right)\right) \leq \max\left\{r_{i-1} + 2a_{i-1}, V\left(j + \frac{l}{N_{i,j}}, x_0\right)\right\} + \frac{1}{2}\min(a_{i-1}, a_i),$$

$$\forall t \in \left[j + \frac{l}{N_{i,j}}, \min(t_{\max}, j+1)\right) \tag{3.58}$$

*Proof of P4:* Let $s \in \{0,1,2,..., N_{i,j} - l - 1\}$ with $j + \frac{l+s}{N_{i,j}} < t_{\max}$. By virtue of (3.57), we distinguish the following two cases:

Case 1: Suppose that

$$V\left(j + \frac{l+s}{N_{i,j}}, x\left(j + \frac{l+s}{N_{i,j}}, j + \frac{l}{N_{i,j}}, x_0; d\right)\right) \geq r_{i-1} \tag{3.59}$$

Then by invoking (3.48b) we get from (3.57), (3.59)

$$0 < V\left(t, x\left(t, j + \frac{l}{N_{i,j}}, x_0; d\right)\right) \leq \max\left\{r_{i-1} + 2a_{i-1}, V\left(j + \frac{l}{N_{i,j}}, x_0\right)\right\} + \frac{1}{2}\min(a_{i-1}, a_i),$$

$$\forall t \in \left[j + \frac{l+s}{N_{i,j}}, j + \frac{l+s+1}{N_{i,j}}\right] \tag{3.60}$$

The desired (3.58) is a consequence of (3.60).

Case 2: Suppose that

$$0 < V\left(j + \frac{l+s}{N_{i,j}}, x\left(j + \frac{l+s}{N_{i,j}}, j + \frac{l}{N_{i,j}}, x_0; d\right)\right) < r_{i-1} \tag{3.61}$$

We show that, when (3.61) holds, then

$$V\left(t, x\left(t, j + \frac{l}{N_{i,j}}, x_0; d\right)\right) \leq r_{i-1} + 2a_{i-1} + \frac{1}{2}\min(a_{i-1}, a_i), \forall t \in \left[j + \frac{l+s}{N_{i,j}}, \min\left(t_{\max}, j + \frac{l+s+1}{N_{i,j}}\right)\right) \tag{3.62}$$

Assume on the contrary that there would exist $t \in \left[j + \frac{l+s}{N_{i,j}}, \min\left(t_{\max}, j + \frac{l+s+1}{N_{i,j}}\right)\right)$ such that

$$V\left(t, x\left(t, j + \frac{l}{N_{i,j}}, x_0; d\right)\right) > r_{i-1} + 2a_{i-1} + \frac{1}{2}\min(a_{i-1}, a_i) \tag{3.63}$$

By (3.61), (3.63), there would exist $t_1 \in \left[j + \frac{l+s}{N_{i,j}}, \min\left(t_{\max}, j + \frac{l+s+1}{N_{i,j}}\right)\right)$ such that

$$V\left(t_1, x\left(t_1, j + \frac{l}{N_{i,j}}, x_0; d\right)\right) = r_{i-1} + 2a_{i-1}; V\left(\xi, x\left(\xi, j + \frac{l}{N_{i,j}}, x_0; d\right)\right) \leq r_{i-1} + 2a_{i-1}, \forall \xi \in \left[j + \frac{l+s}{N_{i,j}}, t_1\right] \tag{3.64}$$

By (3.64) and (3.48b) we get

$$0 < V\left(\xi, x\left(\xi, j+\frac{l}{N_{i,j}}, x_0; d\right)\right) \leq r_{i-1} + 2a_{i-1} + \frac{1}{2}\min(a_{i-1}, a_i), \; \forall \xi \in \left[t_1, \min\left(t_{\max}, j+\frac{l+s+1}{N_{i,j}}\right)\right) \quad (3.65)$$

Combining (3.64), (3.65) we obtain $0 < V\left(\xi, x\left(\xi, j+\frac{l}{N_{i,j}}, x_0; d\right)\right) \leq r_{i-1} + 2a_{i-1} + \frac{1}{2}\min(a_{i-1}, a_i)$ for all $\xi \in \left[j+\frac{l+s}{N_{i,j}}, \min\left(t_{\max}, j+\frac{l+s+1}{N_{i,j}}\right)\right)$, which contradicts hypothesis (3.63).

We conclude from (3.60) and (3.62) that in both cases above we have

$$0 < V\left(t, x\left(t, j+\frac{l}{N_{i,j}}, x_0; d\right)\right) \leq \max\left\{r_{i-1} + 2a_{i-1}, V\left(j+\frac{l}{N_{i,j}}, x_0\right)\right\} + \frac{1}{2}\min(a_{i-1}, a_i),$$

for every $t \in \left[j+\frac{l+s}{N_{i,j}}, \min\left(t_{\max}, j+\frac{l+s+1}{N_{i,j}}\right)\right)$ and for all $s \in \{0,1,2,\ldots, N_{i,j} - l - 1\}$ with $j+\frac{l+s}{N_{i,j}} < t_{\max}$ and the latter implies the desired (3.58). This completes the proof of property P4.

We are now in a position to establish (3.36). Let $(x_0, d) \in (\Re^n \setminus \{0\}) \times M_D$ and $t_0 \in \left[j+\frac{l}{N_{i+1,j}}, j+\frac{l+1}{N_{i+1,j}}\right)$ for some $l \in \{0,1,\ldots, N_{i,j} - 1\}$ with $V(t_0, x_0) \in [r_{i-1}, r_i]$. Then, exploiting inequality (3.48b) we obtain:

$$V(t_0, x_0) - \frac{1}{2}a_i \leq V(t, x(t, t_0, x_0; d)) \leq V(t_0, x_0) + \frac{1}{2}a_i, \; \forall t \in \left[t_0, j+\frac{l+1}{N_{i+1,j}}\right] \quad (3.66)$$

Since $r_i + \frac{1}{2}a_i < r_{i+1} - 2a_{i+1}$, by virtue of (3.66) and (3.59) of P4 we get $0 < V(t, \phi(t, t_0, x_0; d)) \leq r_i + \frac{5}{2}a_i$, for all $t \in [t_0, \min(t_{\max}, j+1))$ and this establishes (3.36). The proof of Lemma 3.2 is complete. ◁

**Lemma 3.3:** *Under the same hypotheses imposed in Lemma 3.1 for system (1.1), there exists a continuous mapping $\tilde{k} : \Re^+ \times (\Re^n \setminus \{0\}) \to U$, being continuously differentiable with respect to $x \in \Re^n \setminus \{0\}$, which satisfies*

$$\left|\tilde{k}(t, x)\right| \leq \tilde{b}(t, x), \; \forall (t, x) \in \Re^+ \times (\Re^n \setminus \{0\}) \quad (3.67)$$

*where $\tilde{b}(\cdot, \cdot)$ is defined by (3.4) and in such a way that that the following property holds for all $(t_0, x_0, d) \in \Re^+ \times (\Re^n \setminus \{0\}) \times M_D$:*

$$V(t, x(t, t_0, x_0; d)) \leq 9V(t_0, x_0), \; \forall t \in [t_0, \min([t_0]+2, t_{\max})) \quad (3.68)$$

*where $x(\cdot, t_0, x_0; d)$ denotes the unique solution of*

$$\dot{x} = f(t, d, x, \tilde{k}(t, x)), (t, x) \in \Re^+ \times (\Re^n \setminus \{0\}) \quad (3.69)$$

with initial condition $x(t_0) = x_0 \in \Re^n \setminus \{0\}$, corresponding to input $d \in M_D$, and $t_{\max} := t_{\max}(t_0, x_0, d) > t_0$ denotes its maximal existence time. Moreover, there exists $\tilde{\rho} \in C^0(\Re^+; \Re^+)$ being positive definite with $\tilde{\rho}(s) \leq s$ for all $s \geq 0$, such that

$$V(2j+2, x(2j+2, 2j, x_0; d)) \leq V(2j, x_0) - \tilde{\rho}(V(2j, x_0)),$$
$$\text{for all } (x_0, d, j) \in (\Re^n \setminus \{0\}) \times M_D \times Z^+ \text{ with } 2j+2 < t_{\max} \tag{3.70}$$

where $t_{\max} > 2j$ in (3.70) is the maximal existence time of the solution $x(\cdot, 2j, x_0; d)$ of (3.69).

**Proof of Lemma 3.3:** Let $r = \{r_i : i \in Z\}$ be a set with $r_i > 0$ and such that

$$r_{i+1} \leq 2r_i \text{ and } \lim_{i \to +\infty} r_i = +\infty, \ \lim_{i \to -\infty} r_i = 0 \tag{3.71}$$

Consider the set

$$r' = \left\{ r'_i = \frac{r_i + r_{i+1}}{2} : i \in Z \right\} \tag{3.72}$$

which by virtue of (3.72) satisfies $r'_i > 0$ and further

$$r'_{i+1} \leq 2r'_i \text{ and } \lim_{i \to +\infty} r'_i = +\infty, \ \lim_{i \to -\infty} r'_i = 0 \tag{3.73}$$

Define

$$\mu_i := \frac{1}{4} \min\{\rho(s) : s \in [r_{i-1}, r_i]\}, \mu'_i := \frac{1}{4} \min\{\rho(s) : s \in [r'_{i-1}, r'_i]\} \tag{3.74}$$

and let $a = \{a_i : i \in Z\}$, $a' = \{a'_i : i \in Z\}$ be a pair of sets satisfying:

$$a_i > 0 \ ; \ a'_i > 0 \tag{3.75a}$$

$$\frac{5}{2} a_i \leq r_{i-1} \ ; \ \frac{5}{2} a'_i \leq r'_{i-1} \tag{3.75b}$$

$$r_i + 2a_i < r_{i+1} - 2a_{i+1} \ ; \ r'_i + 2a'_i < r'_{i+1} - 2a'_{i+1} \tag{3.75c}$$

$$a_i + a'_i \leq \frac{r_{i+1} - r_i}{8} \ ; \ a_i + a'_{i-1} \leq \frac{r_i - r_{i-1}}{8} \tag{3.75d}$$

$$a_i \leq \frac{\mu'_i}{8} \ ; \ a'_i \leq \frac{\mu_{i+1}}{8} \tag{3.75e}$$

By Lemma 3.2, there exist continuous mappings $k_{r,a} : \Re^+ \times (\Re^n \setminus \{0\}) \to U$, $k_{r',a'} : \Re^+ \times (\Re^n \setminus \{0\}) \to U$ being continuously differentiable with respect to $x \in \Re^n \setminus \{0\}$ with

$$k_{r,a}(j, x) = k_{r',a'}(j, x) = 0, \ \forall (x, j) \in (\Re^n \setminus \{0\}) \times Z^+ \tag{3.76a}$$

$$\frac{\partial k_{r,a}}{\partial x}(j, x) = \frac{\partial k_{r',a'}}{\partial x}(j, x) = 0, \ \forall (x, j) \in (\Re^n \setminus \{0\}) \times Z^+ \tag{3.76b}$$

satisfying properties (3.35), (3.36) and (3.38). Finally, consider the map $\tilde{k} : \Re^+ \times (\Re^n \setminus \{0\}) \to U$ defined as:

$$\tilde{k}(t,x) = \begin{cases} k_{r,a}(t,x), & \text{for } t \in [2j, 2j+1), (j,x) \in Z^+ \times (\Re^n \setminus \{0\}) \\ k_{r',a'}(t,x), & \text{for } t \in [2j+1, 2j+2), (j,x) \in Z^+ \times (\Re^n \setminus \{0\}) \end{cases} \quad (3.77)$$

By taking into account (3.76a,b), (3.77) and regularity properties of $k_{r,a}(\cdot,\cdot)$ and $k_{r',a'}(\cdot,\cdot)$, it follows that $\tilde{k}(\cdot,\cdot)$ is continuous, continuously differentiable with respect to $x \in \Re^n \setminus \{0\}$ and satisfies

$$\tilde{k}(j,x) = 0, \ \forall (j,x) \in Z^+ \times (\Re^n \setminus \{0\}) \quad (3.78)$$

Moreover, (3.67) is an immediate consequence of definition (3.77) and inequality (3.35). Also, by (3.36), (3.71) and (3.75b), it follows that for all $(t_0, x_0, d, i) \in \Re^+ \times (\Re^n \setminus \{0\}) \times M_D \times Z$ it holds:

$$V(t_0, x_0) \in [r_{i-1}, r_i] \implies V(t, x_{r,a}(t, t_0, x_0; d)) \le 3V(t_0, x_0), \text{ for all } t \in [t_0, \min([t_0]+1, t_{\max}^{r,a})) \quad (3.79a)$$

$$V(t_0, x_0) \in [r_{i-1}, r_i] \implies V(t, x_{r',a'}(t, t_0, x_0; d)) \le 3V(t_0, x_0), \text{ for all } t \in [t_0, \min([t_0]+1, t_{\max}^{r',a'})) \quad (3.79b)$$

where $x_{r,a}(\cdot, t_0, x_0; d)$ denotes the (unique) solution of

$$\dot{x} = f(t, d, x, k_{r,a}(t,x)), \ (t,x) \in \Re^+ \times (\Re^n \setminus \{0\}) \quad (3.80a)$$

and $x_{r',a'}(\cdot, t_0, x_0; d)$ is the (unique) solution of

$$\dot{x} = f(t, d, x, k_{r',a'}(t,x)) \ (t,x) \in \Re^+ \times (\Re^n \setminus \{0\}) \quad (3.80b)$$

with same initial condition $x(t_0) = x_0 \in \Re^n \setminus \{0\}$ and $d \in M_D$, and $t_{\max}^{r,a} > t_0$ and $t_{\max}^{r',a'} > t_0$, respectively denote their maximal existence times. The desired inequality (3.68) is a direct consequence of (3.79a,b), definition (3.77) and the following obvious fact:

**Fact:** The solution of (3.69) with initial condition $x(t_0) = x_0 \in \Re^n \setminus \{0\}$, corresponding to input $d \in M_D$, is identical for $t \in [t_0, \min([t_0]+1, t_{\max}^{r,a}))$ to the solution $x_{r,a}(t, t_0, x_0; d)$ of (3.80a) if $[t_0]$ is even, and is identical for $t \in [t_0, \min([t_0]+1, t_{\max}^{r',a'}))$ to the solution $x_{r',a'}(t, t_0, x_0; d)$ of (3.80b), if $[t_0]$ is odd.

In order to show (3.70), let $(x_0, d, j) \in (\Re^n \setminus \{0\}) \times M_D \times Z^+$ such that the unique solution $x(\cdot, 2j, x_0; d)$ of (3.69) with initial condition $x(2j) = x_0 \in \Re^n \setminus \{0\}$, corresponding to input $d \in M_D$ is well-defined on $[2j, 2j+2]$ (notice that if there is no such $(x_0, d, j) \in (\Re^n \setminus \{0\}) \times M_D \times Z^+$ then property (3.70) trivially holds for every positive definite function $\tilde{\rho} \in C^0(\Re^+; \Re^+)$). Let $i \in Z$ be the smallest integer with

$$r_{i-1} - 2a_{i-1} < V(2j, x_0) \le r_i - 2a_i \quad (3.81)$$

whose existence is guaranteed from (3.71), (3.75c). By virtue of (3.38), (3.81) and previous fact, it follows that

$$V(2j+1, x(2j+1, 2j, x_0; d)) \le \max(r_{i-1} + 2a_{i-1}, V(2j, x_0) - \mu_i) \quad (3.82)$$

Notice that by virtue of (3.75d), we have $V(2j+1, x(2j+1, 2j, x_0; d)) \le r'_i - 2a'_i$. Consequently, there exists an integer $k \le i$ with

$$r'_{k-1} - 2a'_{k-1} < V(2j+1, x(2j+1, 2j, x_0; d)) \le r'_k - 2a'_k \quad (3.83)$$

We distinguish the following cases:

Case 1: $k < i$

In this case it follows from (3.83) that $V(2j+1, x(2j+1, 2j, x_0; d)) \leq r'_{i-1} - 2a'_{i-1}$. By virtue of (3.38) and the fact above we then obtain

$$V(2j+2, x(2j+2, 2j, x_0; d)) \leq \max\left( r'_{i-2} + 2a'_{i-2}, V(2j, x_0) - \mu'_{i-1} - \mu_i, r_{i-1} + 2a_{i-1} - \mu'_{i-1} \right) \quad (3.84)$$

We now take into account (3.75d) which implies

$$r'_{i-2} + 2a'_{i-2} \leq r_{i-1} - 2a_{i-1} - \frac{r_{i-1} - r_{i-2}}{4} \quad (3.85)$$

From (3.84), (3.85) and the left hand side inequality in (3.81) we get

$$V(2j+2, x(2j+2, 2j, x_0; d)) \leq V(2j, x_0) + \max\left( -\frac{r_{i-1} - r_{i-2}}{4}, 4a_{i-1} - \mu'_{i-1} \right) \quad (3.86)$$

which by virtue of (3.75e) implies:

$$V(2j+2, x(2j+2, 2j, x_0; d)) \leq V(2j, x_0) - \frac{1}{4}\min\left( r_{i-1} - r_{i-2}, 2\mu'_{i-1} \right) \quad (3.87)$$

Case 2: $k = i$.

Notice that, since $r'_{i-1} - 2a'_{i-1} > r_{i-1} + 2a_{i-1}$ (which is a consequence of (3.75d)), we conclude from (3.82) and using the left hand side inequality (3.83) with $k = i$:

$$r'_{i-1} - 2a'_{i-1} + \mu_i < V(2j, x_0) \quad (3.88)$$

Also, by (3.38) and the fact above we get $V(2j+2, x(2j+2, 2j, x_0; d)) \leq \max\left( r'_{i-1} + 2a'_{i-1}, V(2j, x_0) - \mu'_i - \mu_i \right)$, which in conjunction with (3.88) gives $V(2j+2, x(2j+2, 2j, x_0; d)) \leq V(2j, x_0) + 4a'_{i-1} - \mu_i$ and the latter by virtue of (3.75e) implies:

$$V(2j+2, \phi(2j+2, 2j, x_0; d)) \leq V(2j, x_0) - \frac{1}{2}\mu_i \quad (3.89)$$

We conclude from (3.87) and (3.89) that in both cases we have:

$$r_{i-1} - 2a_{i-1} < V(2j, x_0) \leq r_i - 2a_i \Rightarrow V(2j+2, x(2j+2, 2j, x_0; d)) \leq V(2j, x_0) - \gamma_i \quad (3.90a)$$

$$\gamma_i := \frac{1}{4}\min\left( r_{i-1} - r_{i-2}, 2\mu'_{i-1}, 2\mu_i \right) \quad (3.90b)$$

Now let

$$\overline{\rho}(s) := \begin{cases} \frac{\left(\min(\gamma_i, \gamma_{i+1}) - \min(\gamma_{i-1}, \gamma_i)\right)(s - r_{i-1} + 2a_{i-1})}{(r_i - 2a_i - r_{i-1} + 2a_{i-1})} + \min(\gamma_{i-1}, \gamma_i) & , \quad \text{for} \quad s \in (r_{i-1} - 2a_{i-1}, r_i - 2a_i] \\ 0 & , \quad \text{for} \quad s = 0 \end{cases} \quad (3.91)$$

Notice that (3.74), (3.90b) and (3.91) imply that $0 < \min(\gamma_{i-1}, \gamma_i, \gamma_{i+1}) \leq \overline{\rho}(s) \leq \gamma_i$ for $s \in (r_{i-1} - 2a_{i-1}, r_i - 2a_i]$ and further $\lim_{i \to \infty} \mu_i = \lim_{i \to \infty} \mu'_i = \lim_{i \to \infty} \gamma_i = 0$. Thus, we may easily verify that $\overline{\rho} : \Re^+ \to \Re^+$ is positive definite and continuous. Finally, define

$$\tilde{\rho}(s) := \min\{\overline{\rho}(s), s\} \tag{3.92}$$

Property (3.90a,b) in conjunction with (3.91) imply the desired (3.70) is satisfied and the proof is complete. ◁

We are now in a position to prove Theorem 2.8.

**Proof of Theorem 2.8:** By virtue of Lemma 3.3 there exists a continuous mapping $\tilde{k}: \Re^+ \times (\Re^n \setminus \{0\}) \to U$, being continuously differentiable with respect to $x \in \Re^n \setminus \{0\}$, satisfying (3.67), (3.68) and (3.70). We define:

$$K(t, x) = \tilde{k}(t, x) \text{ for } t \geq 0, \ x \neq 0 \tag{3.93a}$$

$$K(t, 0) = 0 \text{ for } t \geq 0 \tag{3.93b}$$

It follows from (2.3), (3.4), (3.67) and definition (3.93) that $K: \Re^+ \times \Re^n \to U$ is a continuous and continuously differentiable mapping with respect to $x \in \Re^n$ on the set $\Re^+ \times (\Re^n \setminus \{0\})$.

**Fact 1:** *For every $(t_0, x_0, d) \in \Re^+ \times \Re^n \times M_D$, the solution $x(\cdot, t_0, x_0; d)$ of (2.6) with initial condition $x(t_0) = x_0$, corresponding to input $d \in M_D$ is unique and is defined for all $t \geq t_0$.*

*Proof of Fact 1:* Consider the resulting system (2.6) with $K(\cdot, \cdot)$ as above and notice that its solution with initial condition $x(t_0) = x_0 \in \Re^n \setminus \{0\}$, corresponding to some $d \in M_D$ coincides with the unique solution of (3.69) evolving on $\Re^+ \times (\Re^n \setminus \{0\})$ with same initial condition $x(t_0) = x_0 \in \Re^n \setminus \{0\}$, and same $d \in M_D$ on the interval $[t_0, t_{\max})$, where $t_{\max} > t_0$ is the maximal existence time of the solution of (3.69). For the case $t_{\max} = +\infty$, the statement of Fact 1 is a direct consequence of previous argument. Suppose next that $t_{\max} < +\infty$. To establish the desired claim, we need the following implication, which is a consequence of (2.1) and (3.68):

$$t_{\max} < +\infty \Rightarrow \lim_{t \to t_{\max}^-} x(t) = 0 \tag{3.94}$$

In order to show (3.94), let $(t_0, x_0, d) \in \Re^+ \times (\Re^n \setminus \{0\}) \times M_D$ and suppose that the maximal existence time $t_{\max} > t_0$ of the (unique) solution of (3.69) with initial condition $x(t_0) = x_0 \in \Re^n \setminus \{0\}$ corresponding to $d \in M_D$ is finite, i.e., $t_{\max} < +\infty$. Repeated use of (3.68) implies that

$$V(t, x(t)) \leq 9^i V(t_0, x_0), \ \forall t \in [t_0, t_{\max})$$

where $i \in Z^+$ is the smallest integer with the property $2i \geq t_{\max}$. The above inequality in conjunction with (2.1) with $\beta(t) \equiv 1$ gives

$$|x(t)| \leq M := \frac{1}{\min_{\tau \in [0, t_{\max}]} \mu(\tau)} a_1^{-1}\left(9^i a_2(|x_0|)\right) < +\infty, \ \forall t \in [t_0, t_{\max}) \tag{3.95}$$

Definition of $t_{\max}$ and (3.95) implies (3.94). By applying standard arguments we may also establish show that for every $(t_0, d) \in \Re^+ \times M_D$, the solution of (2.6) with initial condition $x(t_0) = 0$, corresponding to input $d \in M_D$ is unique and satisfies $x(t) = 0$ for all $t \geq t_0$. Indeed, suppose on the contrary that there exists a nonzero solution of (2.6) with initial condition $x(t_0) = 0$, defined on $[t_0, t_0 + h)$ for some $h > 0$ and let $t_1 \in [t_0, t_0 + h)$ with $x(t_1) \neq 0$ and $a := \max\{t \in [t_0, t_1]: x(t) = 0\}$. Then $x(a) = 0$ and $x(t) \neq 0$ for all $t \in (a, t_1)$. Without loss of generality we may assume that $t_1 \leq [a] + 2$ (if $t_1 > [a] + 2$ then we may use $[a] + 2$ instead of $t_1$, which in this case satisfies $x([a] + 2) \neq 0$). Define $\varepsilon := \frac{1}{2} a_1\left(\mu(t_1) |x(t_1)|\right) > 0$ and let $t_2 \in (a, t_1)$ such that $V(t, x(t)) \leq \frac{\varepsilon}{9}$ for all $t \in [a, t_2]$. By taking into account (3.68) and the fact $t_1 \leq [t_2] + 2$ it then follows that $V(t, x(t)) \leq \varepsilon$ for all $t \in [t_2, t_1]$ and the latter in conjunction with (2.1) yields $a_1(\mu(t_1)|x(t_1)|) \leq \varepsilon$. But this contradicts the definition of $\varepsilon$, hence, we conclude that $x(t) = 0$ for all $t \geq t_0$. The previous discussion in conjunction with (3.94) asserts that the solution $x(\cdot)$ of (2.6) with initial condition $x(t_0) = x_0 \in \Re^n \setminus \{0\}$, corresponding to $d \in M_D$ coincides with the solution of (3.69) with same initial condition, and same $d \in M_D$ on the interval $[t_0, t_{\max})$, $t_{\max} > t_0$ being the maximal existence time of the solution (3.69); moreover, if $t_{\max} < +\infty$, the corresponding solution of (2.6) satisfies $x(t) = 0$ for all $t \geq t_{\max}$ and the proof of Fact 1 is complete.

Fact 1 asserts that, if for some $(t_0, x_0, d) \in \Re^+ \times \Re^n \times M_D$ we consider the maximum existence time $t_{\max} := t_{\max}(t_0, x_0, d)$ of the corresponding solution of (2.6), then $t_{\max} = +\infty$ and the latter in conjunction with (3.68) and (3.70) assert that the following properties are fulfilled for every $(t_0, x_0, d, j) \in \Re^+ \times \Re^n \times M_D \times Z^+$:

$$V(t, x(t, t_0, x_0; d)) \leq 9V(t_0, x_0), \text{ for all } t \in [t_0, [t_0] + 2] \tag{3.96}$$

$$V(2j+2, x(2j+2, 2j, x_0; d)) \leq V(2j, x_0) - \tilde{\rho}(V(2j, x_0)) \tag{3.97}$$

where $\tilde{\rho} \in C^0(\Re^+; \Re^+)$ is the positive definite function involved in (3.70), $x(t, t_0, x_0; d)$ denotes the solution of (2.6) with initial condition $x(t_0) = x_0 \in \Re^n$, corresponding to input $d \in M_D$ and $[t_0]$ is the integer part of $t_0$. The following inequality is a straightforward consequence of inequalities (3.96), (3.97):

$$V(t, x(t, t_0, x_0; d)) \leq 81 V(t_0, x_0), \text{ for all } t \geq t_0 \text{ and } (t_0, x_0, d) \in \Re^+ \times \Re^n \times M_D \tag{3.98}$$

Inequality (3.98) in conjunction with inequality (2.1) with $\beta(t) \equiv 1$, implies Robust Forward Completeness, Uniform Robust Lagrange Output Stability and Uniform Robust Lyapunov Output Stability. Therefore, in order to establish URGAOS for (2.6), it remains to show Uniform Output Attractivity on compact sets of initial states. Let $R \geq 0$, $\varepsilon > 0$ and $(t_0, x_0, d) \in \Re^+ \times \Re^n \times M_D$ with $|x_0| \leq R$ and consider the smallest non-negative integer $j$, which satisfies $t_0 \leq 2j$. Then we have:

**Fact 2:** *For every $\varepsilon > 0$, it holds:*

$$V(2i, x(2i, t_0, x_0; d)) \leq \frac{1}{9} a_1(\varepsilon) \tag{3.99}$$

*for every $i \in Z^+$ with*

$$i \geq j + \frac{9 a_2(R)}{r}; \; r := \min\left\{\tilde{\rho}(s) : s \in \left[\frac{1}{9} a_1(\varepsilon), \frac{1}{9} a_1(\varepsilon) + 9 a_2(R)\right]\right\} \tag{3.100}$$

*Proof of Fact 2:* Suppose on the contrary that there exists $\varepsilon > 0$, $i \geq j + \frac{9 a_2(R)}{r}$ with $V(2i, x(2i, t_0, x_0; d)) > \frac{1}{9} a_1(\varepsilon)$. Using (3.97), it follows that $V(2k, x(2k, t_0, x_0; d)) > \frac{1}{9} a_1(\varepsilon)$ for $k = j, ..., i$. Also (3.97) implies that $V(2k, x(2k, t_0, x_0; d)) \leq V(2j, x(2j, t_0, x_0; d))$, for $k = j, ..., i$. Consequently, from (2.1) (with $\beta(t) \equiv 1$), (3.96) we get

$$V(2k, x(2k, t_0, x_0; d)) \in \left[\frac{1}{9} a_1(\varepsilon), \frac{1}{9} a_1(\varepsilon) + 9 a_2(R)\right], \text{ for } k = j, ..., i \tag{3.101}$$

On the other hand, by recalling (3.97) and using (3.100), (3.101) it follows:

$$V(2(k+1), x(2(k+1), t_0, x_0; d)) \leq V(2k, x(2k, t_0, x_0; d)) - r, \text{ for } k = j, ..., i$$

which in turns gives:

$$V(2k, x(2k, t_0, x_0; d)) \leq V(2j, x(2j, t_0, x_0; d)) - r(k-j), \text{ for } k = j,...,i \tag{3.102}$$

Using (3.96) and (2.1) with $\beta(t) \equiv 1$ we get:

$$V(2j, x(2j, t_0, x_0; d)) \leq 9a_2(R) \tag{3.103}$$

Inequalities (3.102), (3.103) in conjunction with the fact that $i \geq j + \dfrac{9a_2(R)}{r}$, give $V(2i, x(2i, t_0, x_0; d)) \leq 0$, which contradicts the hypothesis $V(2i, x(2i, t_0, x_0; d)) > \dfrac{1}{9} a_1(\varepsilon)$ and the proof of Fact 2 is complete.

Applying again (2.1) with $\beta(t) \equiv 1$ and (3.96), (3.99) of Fact 2, it follows that for every $R \geq 0$, $\varepsilon > 0$, $(t_0, x_0, d) \in \Re^+ \times \Re^n \times M_D$ with $|x_0| \leq R$, it holds that:

$$|H(t, x(t, t_0, x_0; d))| \leq \varepsilon \text{ for all } t \geq t_0 + 2 + \frac{18 a_2(R)}{r}$$

where $r := \min\left\{ \tilde{\rho}(s) : s \in \left[ \dfrac{1}{9} a_1(\varepsilon), \dfrac{1}{9} a_1(\varepsilon) + 9a_2(R) \right] \right\}$ and this establishes Uniform Output Attractivity on compact sets of initial states. The proof is complete. ◁

For the proof of Theorem 2.9 we need an additional lemma, which provides sufficient conditions for (non-uniform in time) RGAOS. It is important to mention here the paper [17], where, under different hypotheses than those imposed below, asymptotic stability for time-varying system is explored by estimating the difference between values of an appropriate Lyapunov function along the trajectories of system at a given sequence of times.

**Lemma 3.4:** *Consider system (2.4), where $f : \Re^+ \times D \times \Re^n \to \Re^n$, $H : \Re^+ \times \Re^n \to \Re^k$ are continuous and for every bounded interval $I \subset \Re^+$ and every compact set $S \subset \Re^n$ there exists $L \geq 0$ such that $|f(t, d, x) - f(t, d, y)| \leq L|x - y|$ for all $(t, d) \in I \times D$, $x, y \in S$. Moreover, assume that the set $D \subset \Re^l$ is compact and $f(t, d, 0, 0) = 0$, $H(t, 0) = 0$ for all $(t, d) \in \Re^+ \times D$. Suppose that there exists a function $\gamma \in K^+$ satisfying*

$$\sum_{j=0}^{+\infty} \gamma(2j) < +\infty \tag{3.104a}$$

$$\lim_{t \to +\infty} \gamma(t) = 0 \tag{3.104b}$$

*and further there exist functions $V \in C^1(\Re^+ \times \Re^n; \Re^+)$, $a_1, a_2, a \in K_\infty$, $\mu, \beta, \gamma \in K^+$, $\rho \in C^0(\Re^+; \Re^+)$ being positive definite such that (2.1) holds and the following properties are fulfilled for all $(t_0, x_0, d, j) \in \Re^+ \times \Re^n \times M_D \times Z^+$:*

$$\sup_{t \in [t_0, [t_0]+2]} V(t, x(t, t_0, x_0; d)) \leq a(V(t_0, x_0)) + \gamma(t_0) \tag{3.105a}$$

$$V(2j+2, x(2j+2, 2j, x_0; d)) \leq V(2j, x_0) - \rho(V(2j, x_0)) + \gamma(2j) \tag{3.105b}$$

*where $x(\cdot, t_0, x_0; d)$ denotes the unique solution of (2.4) with initial condition $x(t_0) = x_0 \in \Re^n$, corresponding to input $d \in M_D$. Then system (2.4) is RGAOS.*

**Proof of Lemma 3.4:** Let $(t_0, x_0, d) \in \Re^+ \times \Re^n \times M_D$ and let $j \in Z^+$ be the smallest integer, which satisfies $t_0 \leq 2j$. Inequality (3.105b) implies that

$$V(2i, x(2i, t_0, x_0, d)) \leq V(2j, x(2j, t_0, x_0, d)) + \sum_{k=0}^{i} \gamma(2k), \text{ for all integers } i \geq j \qquad (3.106)$$

Let $M := \sum_{k=0}^{+\infty} \gamma(2k)$ and $B := \sup_{t \geq 0} \gamma(t)$. Then by (2.1), (3.105a) and (3.106) we get:

$$V(t, x(t, t_0, x_0, d)) \leq a\big(a\big(a_2(\beta(t_0)|x_0|)\big) + B + M\big) + B, \text{ for all } t \geq t_0 \qquad (3.107)$$

Inequality (3.107) in conjunction with (2.1) implies RFC and Robust Lagrange Output Stability. Therefore, according to Lemma 3.5 in [12], in order to establish RGAOS, it suffices to show that system (2.4) satisfies the property of Uniform Output Attractivity on compact sets of initial data. To establish this property, consider arbitrary constants $\varepsilon > 0$, $R \geq 0$, $T \geq 0$ and let $(t_0, x_0, d) \in \Re^+ \times \Re^n \times M_D$ with $t_0 \in [0, T]$ and $|x_0| \leq R$. Define $K := a\big(a\big(a_2(R \max_{t \in [0,T]} \beta(t))\big) + B + M\big) + B$. Then by (3.107) it holds that

$$V(t, x(t, t_0, x_0, d)) \leq K, \quad \forall t \geq t_0 \qquad (3.108)$$

Define

$$\tilde{\rho}(s) := \min_{s \leq y \leq K} \rho(y) \qquad (3.109)$$

which obviously is a non-decreasing and continuous function and let $J \geq 0$ be an integer with $\frac{1}{2}\tilde{\rho}(K) \geq \gamma(2i)$ for all integers $i \geq J$, whose existence is guaranteed from (3.104b). Define the sequence

$$q_i := \inf\left\{s \in [0, K] : \frac{1}{2}\tilde{\rho}(s) \geq \gamma(2i)\right\} \text{ for } i \geq J \qquad (3.110)$$

Notice, by virtue of (3.104b) and (3.110) that $q_i \to 0$ and consequently, there exists an integer $N := N(\varepsilon, K) \geq J$ such that

$$q_i + \gamma(2i) \leq S(\varepsilon) \text{ and } \gamma(2i) \leq \frac{1}{2}a_1(\varepsilon), \text{ for all } i \geq N \qquad (3.111)$$

where $a_1 \in K_\infty$ is the function involved in (2.1) and $S(\varepsilon) > 0$ is defined by

$$S(\varepsilon) := a^{-1}\left(\frac{1}{2}a_1(\varepsilon)\right) \qquad (3.112)$$

Notice next that (3.105b) asserts that for all integers $i \geq \max(N, j)$ the following holds:

$$V(2(i+1), x(2(i+1), t_0, x_0; d)) \leq \max\left(S(\varepsilon), V(2i, x(2i, t_0, x_0; d)) - \frac{1}{2}\rho(V(2i, x(2i, t_0, x_0; d)))\right) \qquad (3.113)$$

Indeed, to establish (3.113) we may distinguish two cases. First assume that $V(2i, x(2i, t_0, x_0, d)) \geq q_i$. Then it follows from (3.108), (3.109) and (3.110) that $\frac{1}{2}\rho(V(2i, x(2i, t_0, x_0, d))) \geq \gamma(2i)$ and the latter in conjunction with (3.105b) implies (3.113). The other case is $V(2i, x(2i, t_0, x_0, d)) \leq q_i$. Then the latter in conjunction with (3.105b) and (3.111) implies again (3.113).

The following is a consequence of (3.113):

$$V(2i, x(2i, t_0, x_0, d)) \leq S(\varepsilon) \text{ for all integers } i \geq \max(N, j) + \frac{2K}{\tilde{\rho}(S(\varepsilon))} + 1 \quad (3.114)$$

To show (3.114), suppose on the contrary that there exists integer $i \geq \max(N, j) + \frac{2K}{\tilde{\rho}(S(\varepsilon))} + 1$ with $V(2i, x(2i, t_0, x_0, d)) > S(\varepsilon)$. Then (3.113) implies that

$$V(2k, x(2k, t_0, x_0, d)) > S(\varepsilon) \text{ for all } k = \max(N, j), \max(N, j) + 1, \ldots, i \quad (3.115)$$

From (3.108), (3.109), (3.113), (3.115) it follows that

$$V(2(k+1), x(2(k+1), t_0, x_0, d)) \leq V(2k, x(2k, t_0, x_0, d)) - \frac{1}{2} \tilde{\rho}(S(\varepsilon)), \text{ for all } k = \max(N, j), \max(N, j) + 1, \ldots, i-1$$

which directly implies

$V(2k, x(2k, t_0, x_0, d)) \leq K - (k - \max(N(\varepsilon, K), j)) \frac{1}{2} \tilde{\rho}(S(\varepsilon))$, for all $k = \max(N, j), \max(N, j) + 1, \ldots, i$. The previous inequality for $k = i$ gives $V(2i, x(2i, t_0, x_0, d)) < 0$, which is a contradiction, hence (3.114) is established.

Using (3.105a) and (3.114) we obtain $\sup_{t \in [2i, 2i+2]} V(t, \phi(t, t_0, x_0; d)) \leq a(S(\varepsilon)) + \gamma(2i)$, for all integers $i \geq \max(N, j) + \frac{2K}{\tilde{\rho}(S(\varepsilon))} + 1$. This in conjunction with (3.111) and (3.112) gives:

$$\sup_{t \geq 2i} V(t, \phi(t, t_0, x_0; d)) \leq a_1(\varepsilon), \text{ for all integers } i \geq \max(N, j) + \frac{2K}{\tilde{\rho}(S(\varepsilon))} + 1$$

Using the inequality above and (2.1), we may conclude that the property of Uniform Output Attractivity on compact sets of initial data holds for system (2.4). This completes the proof of Lemma 3.4. ◁

We are now in a position to prove Theorem 2.9.

**Proof of Theorem 2.9:** According to the statement of Lemma 3.3 there exists a continuous mapping $\tilde{k}: \Re^+ \times (\Re^n \setminus \{0\}) \to U$, being continuously differentiable with respect to $x \in \Re^n \setminus \{0\}$, which satisfies (3.67), (3.68), (3.70). We define:

$$K(t, x) := h\left(\frac{V(t, x) - \exp(-t)}{\exp(-t)}\right) \tilde{k}(t, x), \text{ for } V(t, x) > \exp(-t) \quad (3.116a)$$

$$K(t, x) := 0, \text{ for } V(t, x) \leq \exp(-t) \quad (3.116b)$$

where $h: \Re \to [0, 1]$ is a smooth non-decreasing function with $h(s) = 0$ for $s \leq 0$ and $h(s) = 1$ for $s \geq 1$. It can be easily verified that, according to definition (3.116) and the properties of $\tilde{k}: \Re^+ \times (\Re^n \setminus \{0\}) \to U$, the map $K$ takes values in $U$ and satisfies $K(t, 0) = 0$ for all $t \geq 0$. Moreover, $K: \Re^+ \times \Re^n \to U$, is a continuous and continuously differentiable mapping with respect to $x \in \Re^n$ on $\Re^+ \times \Re^n$.

In order to prove Theorem 2.9 we will make use of Lemma 3.4 and three facts below concerning certain properties of the solution of (2.6). Let

$$T = T(t_0, x_0, d) := \begin{cases} \inf\{t \geq t_0 : \exp(t)V(t, x(t)) < 2\} & \text{if} \quad \{t \geq t_0 : \exp(t)V(t, x(t)) < 2\} \neq \emptyset \\ +\infty & \text{if} \quad \{t \geq t_0 : \exp(t)V(t, x(t)) < 2\} = \emptyset \end{cases} \quad (3.117)$$

where $x(\cdot) = x(\cdot, t_0, x_0; d)$ denotes the unique solution of (2.6) with initial condition $x(t_0) = x_0 \in \Re^n$ corresponding to some $d \in M_D$. The following fact is an immediate consequence of (3.116), (3.117) and continuity of the mapping $t \to V(t, x(t))$.

**Fact 1:** *The unique solution $x(\cdot) = x(\cdot, t_0, x_0; d)$ of (2.6) with initial condition $x(t_0) = x_0 \in \Re^n \setminus \{0\}$, satisfying $V(t_0, x_0) \geq 2\exp(-t_0)$, corresponding to some $d \in M_D$ coincides with the unique solution of (3.69) with same initial condition and same $d \in M_D$ on the interval $[t_0, T]$, where $T = T(t_0, x_0, d)$ is defined by (3.117) and*

$$V(T, x(T)) = 2\exp(-T) \text{ if } \{t \geq t_0 : \exp(t)V(t, x(t)) < 2\} \neq \emptyset \quad (3.118)$$

Next, we prove the following:

**Fact 2:** *For the system (2.6), the following property holds for all $(j, x_0, d) \in Z^+ \times \Re^n \times M_D$:*

$$V(2j+2, x(2j+2, 2j, x_0; d)) \leq V(2j, x_0) - \tilde{\rho}(V(2j, x_0)) + 18\exp(-2j) \quad (3.119)$$

*Proof of Fact 2:* Obviously, the desired (3.119) holds for $x_0 = 0$. Next, assume that $x_0 \neq 0$. Let $t_{\max} > 2j$ the maximal existence time of $x(\cdot, 2j, x_0; d)$. We distinguish two cases. The first case is

$$\{t \in [2j, \min(t_{\max}, 2j+2)) : \exp(t)V(t, x(t, 2j, x_0; d)) < 2\} = \emptyset \quad (3.120)$$

In this case, Fact 1 in conjunction with inequalities (2.1), (3.68) and (3.70) guarantee that $t_{\max} > 2j+2$ and that (3.119) holds. The second case is

$$\{t \in [2j, \min(t_{\max}, 2j+2)) : \exp(t)V(t, x(t, 2j, x_0; d)) < 2\} \neq \emptyset \quad (3.121)$$

Let

$$t_1 := \sup\{t \in [2j, \min(t_{\max}, 2j+2)) : \exp(t)V(t, x(t, 2j, x_0; d)) < 2\} \quad (3.122)$$

Clearly, we have from (3.122)

$$\limsup_{t \to t_1^-} \exp(t)V(t, x(t, 2j, x_0; d)) \leq 2 \quad (3.123)$$

and this by virtue of (2.1) implies $t_{\max} > t_1$. If $t_1 = 2j+2$ the desired (3.119) follows from (3.123) holds. If $t_1 < 2j+2$, definition (3.122) guarantees that $\exp(t)V(t, x(t, 2j, x_0; d)) \geq 2$ for all $t \in [t_1, \min(t_{\max}, 2j+2))$ and the latter in conjunction with (3.123) gives

$$\exp(t_1)V(t_1, x(t_1, 2j, x_0; d)) = 2 \quad (3.124)$$

Using Fact 1 together with (3.68) and (3.124) we get $V(t, x(t, 2j, x_0; d)) \leq 18\exp(-t_1)$ for all $t \in [t_1, \min(t_{\max}, 2j+2))$. By exploiting (2.1) we conclude that $t_{\max} > 2j+2$ and therefore the estimate $V(t, x(t, 2j, x_0; d)) \leq 18\exp(-t_1)$ is fulfilled for every $t \in [t_1, 2j+2]$. The latter implies (3.119) and this completes proof of Fact 2.

Finally we show the following fact.

**Fact 3:** *The following property holds for system (2.6):*

$$V(t, x(t, t_0, x_0; d)) \leq 9V(t_0, x_0) + 18\exp(-t_0), \text{ for all } t \in [t_0, [t_0]+2], (t_0, x_0, d) \in \Re^+ \times \Re^n \times M_D \quad (3.125)$$

*Proof of Fact 3:* Obviously, (3.125) holds for $x_0 = 0$. Suppose next that $x_0 \neq 0$ and let us on the contrary assume that there exists $\hat{t} \in [t_0, [t_0]+2]$ with

$$V(\hat{t}, x(\hat{t}, t_0, x_0; d)) > 9V(t_0, x_0) + 18\exp(-t_0) \quad (3.126)$$

We distinguish two cases. First assume that

$$\{s \in [t_0, \hat{t}] : \exp(s)V(s, x(s, t_0, x_0; d)) < 2\} = \emptyset$$

In this case, (3.68) guarantees that $V(\hat{t}, x(\hat{t}, t_0, x_0; d)) \leq 9V(t_0, x_0)$, which contradicts (3.126). Consider the remaining case

$$\{s \in [t_0, \hat{t}] : \exp(s)V(s, x(s, t_0, x_0; d)) < 2\} \neq \emptyset$$

and let $t_1 := \sup\{s \in [t_0, \hat{t}] : \exp(s)V(s, x(s, t_0, x_0; d)) < 2\}$. If $t_1 = \hat{t}$, we would have $V(\hat{t}, x(\hat{t}, t_0, x_0; d)) \leq 2\exp(-t_0)$, which contradicts (3.126). If $t_1 < \hat{t}$, then we would have $\exp(s)V(s, x(s, t_0, x_0; d)) \geq 2$ for all $s \in [t_1, \hat{t}]$ and $\exp(t_1)V(t_1, x(t_1, 2j, x_0; d)) = 2$. Therefore (3.68) gives $V(s, \phi(s, t_0, x_0; d)) \leq 18\exp(-t_1)$ for all $s \in [t_1, \hat{t}]$, which again contradicts (3.126) and we conclude that (3.125) holds. This completes proof of Fact 3.

Inequalities (3.119), (3.125) in conjunction with Lemma 3.4 show that (2.6) is RGAOS and the proof of Theorem 2.9 is complete. ◁

## 4. Conclusions

For general time-varying systems, it is established that existence of an "Output Robust Control Lyapunov Function" implies existence of continuous time-varying feedback stabilizer, which guarantees output asymptotic stability with respect to the resulting closed-loop system. The main results of the present work constitute generalizations of a well known result towards feedback stabilization due to J. M. Coron and L. Rosier in [7] concerning stabilization of autonomous systems by means of time-varying periodic feedback. Further extensions towards same subject, including stabilization of time-varying systems (1.1) by means of discontinuous time-varying feedback in the Fillipov sense (see [3,8,24]) and existence of appropriate control Lyapunov functions will be a subject of forthcoming research.